\documentclass[12pt]{article}

\usepackage {amsmath, amssymb, amsthm}

\usepackage {tikz}
\usetikzlibrary{calc}

\usepackage{hyperref}
\hypersetup{
    colorlinks = true,
    citecolor = red,
}

\usepackage [margin=1.5in] {geometry}

\theoremstyle {definition}

\newtheorem {ex} {Example}

\theoremstyle {plain}
\newtheorem {lemma}{Lemma}
\newtheorem {thm}{Theorem}
\newtheorem {prop}[thm]{Proposition}
\newtheorem {conj}{Conjecture}

\newcommand{\C}{{\mathcal C}}

\newcommand{\G}{{\mathcal G}}
\newcommand{\be}{\begin{equation}}
\newcommand{\ee}{\end{equation}} 

\newcommand{\eps}{\epsilon}

\newcommand{\R}{{\mathbb R}}
\newcommand{\Z}{{\mathbb Z}}
\newcommand{\T}{{\mathbb T}}

\newcommand{\SL}{\mathrm{SL}}
\newcommand{\GL}{\mathrm{GL}}
\newcommand{\Gr}{\mathrm{Gr}}
\newcommand{\Grnn}{\mathrm{Gr}^{\ge}}
\newcommand{\Grtp}{\mathrm{Gr}^{+}}
\newcommand{\bW}{W_{\partial}}
\newcommand{\bB}{B_{\partial}}
\newcommand{\Tr}{\mathrm{Tr}}
\newcommand{\tr}{\mathrm{tr}}
\newcommand{\n}{\mathbf{n}}
\renewcommand{\H}{\mathcal{H}}

\newcommand{\old}[1]{}

\begin {document}
\title {Higher-Rank Dimer Models}
\author{Richard Kenyon\footnote{Department of Mathematics, Yale University, New Haven; richard.kenyon at yale.edu} \and
Nicholas Ovenhouse\footnote{Department of Mathematics, Yale University, New Haven}}
\date{}
\maketitle

\begin{abstract}
Let $\G$ be a bipartite planar graph with edges directed from black to white.
For each vertex $v$ let $n_v$ be a positive integer. A \emph{multiweb} in $\G$ is a multigraph with multiplicity
$n_v$ at vertex $v$. A \emph{connection} is a choice of linear maps on edges $\Phi=\{\phi_{bw}\}_{bw\in E}$
where $\phi_{bw}\in \mathrm{Hom}(\R^{n_b},\R^{n_w})$.  Associated to $\Phi$ is a function on multiwebs, the trace $\Tr_{\Phi}$.
We define an associated Kasteleyn matrix $K=K(\Phi)$ in this setting and write $\det K$ as the sum of traces of all multiwebs. 
This generalizes Kasteleyn's theorem and the result of \cite{DKS}. 

We study connections with positive traces, and define the associated probability measure on multiwebs.
By careful choice of connection we can thus encode the ``free fermionic'' subvarieties
for vertex models such as the $6$-vertex model and $20$-vertex models, and in particular give
determinantal solutions.

We also find for each multiweb system an equivalent scalar system, that is, a planar bipartite graph $H$ and a local measure-preserving mapping
from dimer covers of $H$ to multiwebs on $\G$. We identify a family of positive connections
as those whose scalar versions have positive face weights. 
\end{abstract}

\tableofcontents

\section {Introduction}

The \emph{dimer model} is the study of random dimer covers (also known as perfect matchings) of graphs.
Kasteleyn showed in \cite{Kast} that one can enumerate (weighted) dimer covers of a planar graph
with the determinant of an associated matrix, now called the \emph{Kasteleyn matrix}. The dimer model,
originally defined as a stat mech model, has turned out to have connections
with many other parts of mathematics, including cluster algebras \cite{speyer, ms_10}, representation theory \cite{DKS, KenyonShi}, geometry \cite{KOS, KLRR}, integrable systems \cite{gk_13},
and string theory \cite{Hananyetal}.

The recent works \cite{DKS} and \cite{KenyonShi} discuss relations between the dimer model and $\GL_n$-connections on graphs. 
In particular associated to a bipartite planar graph $\G$ with a $\GL_n$-connection is a certain matrix,
the $\GL_n$-Kasteleyn matrix, whose determinant counts traces of multiwebs in $\G$ (see definitions below).
We generalize here the result of \cite{DKS} from the case of $\GL_n$-connections 
to the case of a ``quiver representation": a graph of vector spaces and linear maps
between them. See Theorem \ref{main} below.

One application of this result is for planar, biperiodic graphs, with bi-periodic connections. 
We give conditions (see Section \ref{scalarsection}) under which all traces of multiwebs are nonnegative, so that one can discuss Gibbs measures on 
random multiwebs. See for example Theorem \ref{p2} for a result regarding random $2$-webs on the honeycomb
with a periodic $\GL_2$-connection.
Using the Postnikov/Talaska parameterization of the totally nonnegative 
Grassmannian $\Grnn_{k,n}$, see \cite{postnikov_06, talaska},
we also construct an explicit, local, measure-preserving mapping between 
multiwebs on $\G$ and dimer covers of a related planar graph.   See Theorem \ref{measurepres} and the example in Figures \ref{hatG2} and \ref{Gdimerex}.

As another application we are able to parameterize and solve (in the sense of computing the free energy,
local statistics, limit shapes, etc.)
the so-called free-fermionic subvarieties of vertex models such as the $6$-vertex and $20$-vertex models, 
in a conceptually simple manner. See Sections \ref{6vsection} and \ref{20vsection}.
\begin{figure}
\begin{center}
\includegraphics[width=1.5in]{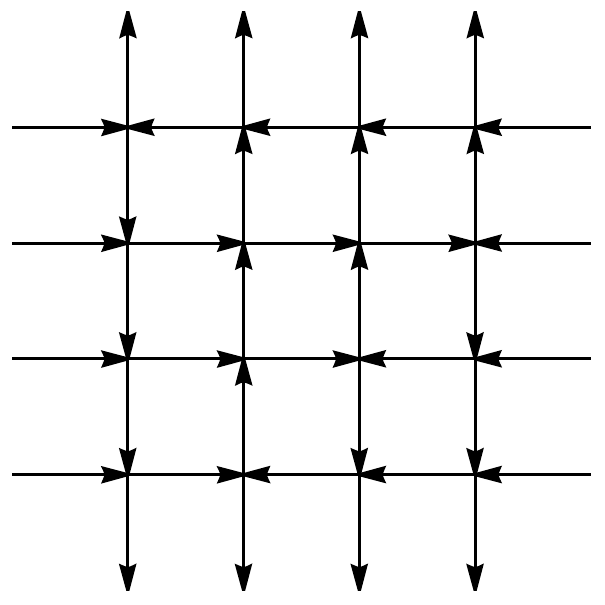}\hskip1cm\includegraphics[width=1.5in]{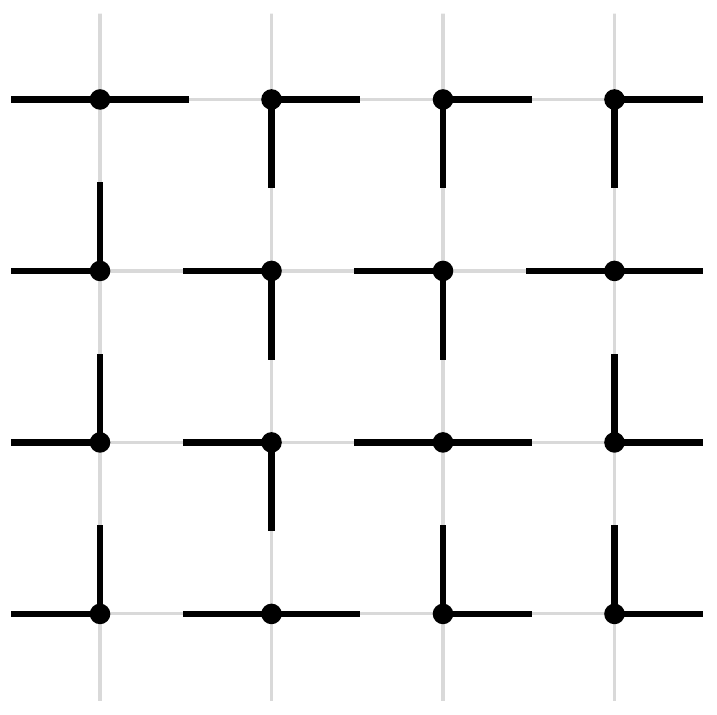}
\end{center}
\caption{\label{6V}Left: orientation version of $6$-vertex model. Right: square ice version.
The bijection is given by the proximity of the arrowheads to the vertices.}
\end{figure}
\old{
\begin{figure}
\begin{center}\includegraphics[width=1.5in]{6Vexample.pdf} 
\end{center}
\caption{\label{6Vexample}A configuration in the $6$-vertex model.}
\end{figure}
}

\begin{figure}
\begin{center}\includegraphics[width=3.5in]{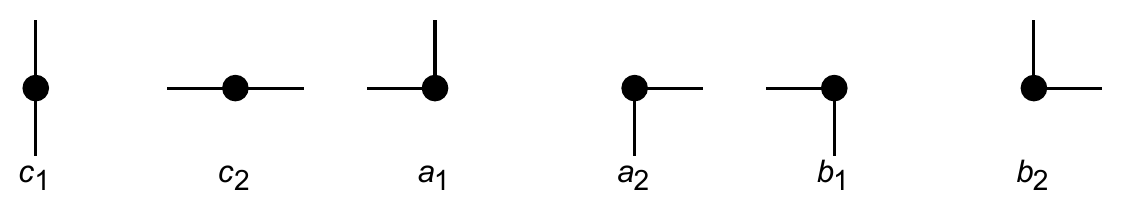} 
\end{center}
\caption{\label{6Vweights}Local weights for the $6$-vertex model.}
\end{figure}

Recall that the six-vertex model is a probability measure on orientations of $\Z^2$ (or subgraphs of $\Z^2$)
with the constraint that at each vertex two edges point in and two point out; see Figure \ref{6V}, left. 
Equivalently, one can consider ``square ice" configurations in which water molecules
are arranged so that an oxygen atom sits at each point of $\Z^2$ and a hydrogen atom on each edge,
so that each oxygen is bonded to two of the neighboring hydrogens (Figure \ref{6V}, right).
In either version there are $6$ possible local configurations at a vertex. By definition, the probability of a configuration is proportional to the product of local weights $c_1,c_2,a_1,a_2,b_1,b_2$
at each vertex which depend on the structure there; see Figure \ref{6Vweights}.
Computing the free energy in the thermodynamic limit (that is, on $\Z^2$) is a long standing open problem.
Lieb \cite{Lieb} computed the free energy in the symmetric case $a_1=a_2,b_1=b_2,c_1=c_2$.
Another subvariety of the space of vertex weights, the $5$-vertex model (where one of $a_1,a_2,b_1,b_2$ is set to $0$) was recently solved in \cite{dGKW}. 
A third subcase, the free-fermionic subvariety, is defined by the constraint that $a_1a_2+b_1b_2-c_1c_2=0$. 
This case is solvable by determinantal methods, see e.g. \cite{Baxter}. 
This last case, and the more general
free fermionic 6V models with staggered weights (that is, invariant under a coarser sublattice), can be solved with our methods.

\begin{figure}
\begin{center}\includegraphics[width=1.5in]{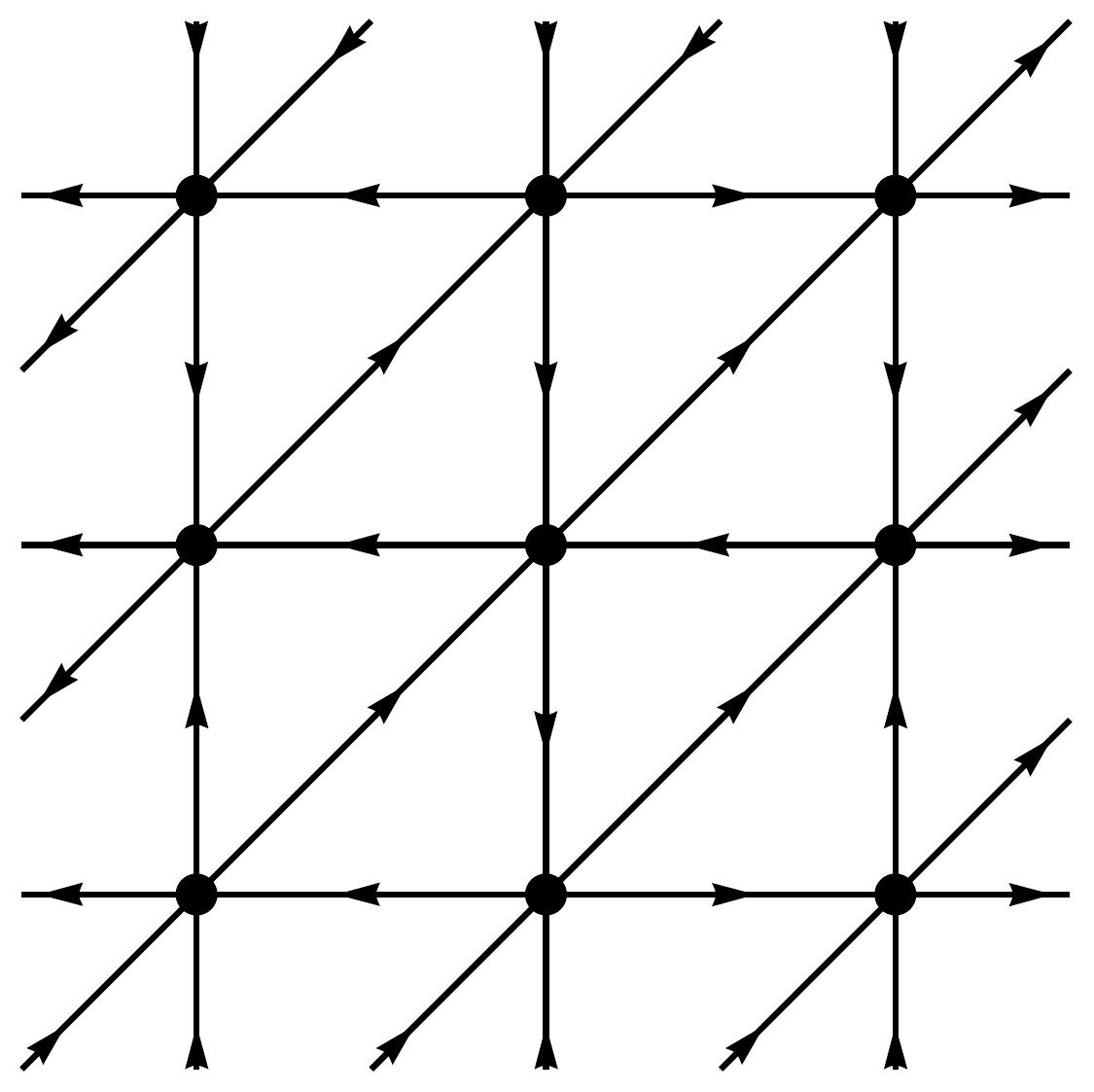}\hskip1cm \includegraphics[width=1.5in]{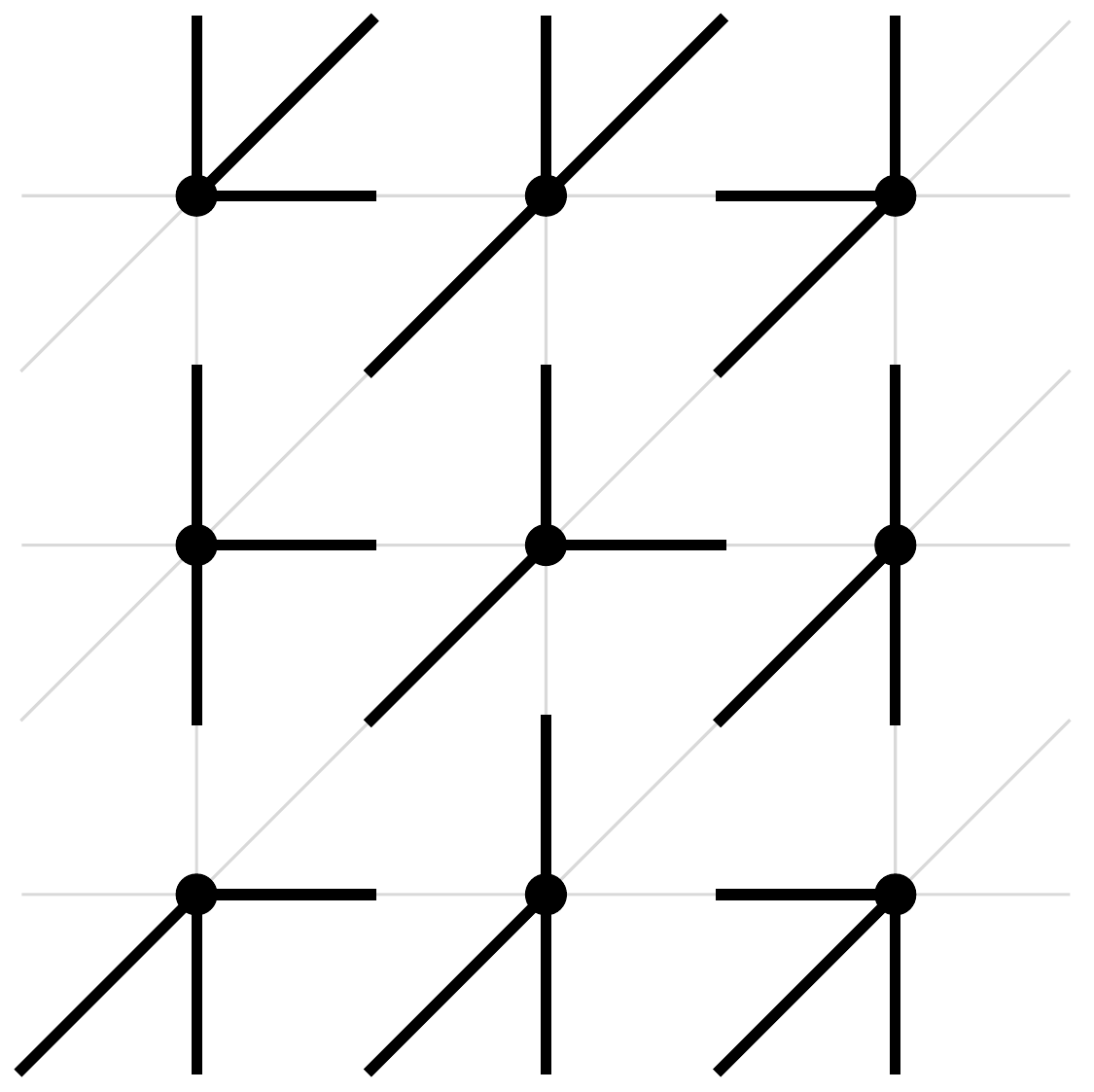} 
\end{center}
\caption{\label{20Vexample}A configuration in the $20$-vertex model.}
\end{figure}

In a similar manner the $20$-vertex model (Figure \ref{20Vexample})
 is a model of orientations of the triangular grid (i.e. the square grid with diagonals) with indegree $3$ at each vertex; it is defined by $20$ local vertex weights.
It also has a free-fermionic subvariety, which is a certain subvariety of dimension $9$, discussed below. 
Again in this case we show that the 
model restricted to this subvariety is solvable by determinantal methods. 
Vertex models of other multiplicities and on other planar graphs can likewise be solved by similar means. 

As a final application, the dimer model on a bipartite graph on a torus has an associated Liouville integrable
system \cite{gk_13}. Pulling this system back from the scalarized version to the quiver representation gives a novel integrable structure on
such quiver representations. We intend to explore this connection in a follow-up paper \cite{KO2}.
\bigskip

\noindent{\bf Acknowledgments.} This research was supported by NSF grant DMS-1940932 
and the Simons Foundation grant 327929. 
We thank Daniel Douglas, Sam Panitch and Sri Tata for discussions.

\section {Webs and Tensor Networks}

In this section, we review and generalize the definitions of tensor networks, webs, multiwebs and multiweb traces. 
Background can be found in \cite{fp_16, kuperberg, DKS}.

For a positive integer $n$, a \emph{planar $n$-web} is a planar bipartite $n$-regular graph.
Such webs are used in the study of the representation theory of $\SL_n$ \cite{sikora_01}.
The definition of web can be generalized to work inside a given graph \cite{DKS}, where it is called
a \emph{multiweb}. We can also allow $n$ to vary from vertex to vertex (Postnikov \cite{postnikov_18} also discusses this generalization from another viewpoint).

\subsection{$\n$-Multiwebs}

Let $\G=(V=B\cup W,E)$ be a bipartite graph.
Let $\n:V\to\Z_{>0}$ be a positive integer function on vertices, called the vertex multiplicity.
An \emph{$\n$-multiweb in $\G$} is a function $m:E\to\Z_{\ge 0}$ that sums to $n_v=\n(v)$ 
at each vertex $v$, that is, 
such that $\sum_{u\sim v} m_{uv} = n_v$.
We call $m_{e}$ the \emph{multiplicity} of the edge $e$. The multiplicity of an edge may be zero.
Let $\Omega_{\n}=\Omega_{\n}(\G)$ be the set of $\n$-multiwebs. 
When $n_v\equiv 1$, a $1$-multiweb is also known as a \emph{dimer cover}, or \emph{perfect matching}, of $\G$.
When $n_v\equiv 2$, a $2$-multiweb is also known as a \emph{double-dimer cover}, see for example \cite{KenyonWilson}.

The set $\Omega_{\n}$ is nonempty only if $\n$ satisfies 
$\sum_{b\in B}n_b = \sum_{w\in W}n_w$ (by the handshake lemma) and also 
certain linear inequalities. These inequalities are given by Hall's matching theorem:

\begin{thm}[\cite{hall_35}] 
    $\Omega_{\n}$ is nonempty if and only if $\sum_{b\in B}n_b = \sum_{w\in W}n_w$ and for every set $X$ of black vertices, there are sufficiently many white neighbors, 
    counted with multiplicity: for any set $X\subset B$, we have
    $\sum_{b\in X}n_b \le \sum_{w\in Y}n_w$ where $Y$ is the set of neighbors of $X$.
\end{thm}

\begin{proof} 
    Replace each vertex $v$ with $n_v$ vertices and each edge $bw$ with a complete bipartite graph
    $K_{n_b,n_w}$. Then the statement follows immediately from Hall's matching theorem \cite{hall_35} regarding the existence of a perfect matching of a bipartite graph.
\end{proof}


\subsection{Graphs with boundary}

Even though our main theorem, Theorem \ref{main}, deals with graphs without boundary,
in principle there is a generalization to graphs with boundary (with a slightly more technical statement and proof which we do not give here); see \cite{KenyonShi} which deals 
with the case of graphs with boundary in the case $\n\equiv 3$. 

Suppose $\G$ is planar, embedded in the disk, with a distinguished set of vertices $V_\partial=\bW\cup \bB$ on its outer face,
called boundary vertices.
Fix $\n:V\to\Z_{>0}$ as before.
Also fix a ``boundary multiplicity function" $d:V_\partial\to\Z_{>0}$ satisfying $d_v\le n_v$. 
We define an $(\n,d)$-multiweb as above except that at a boundary vertex $v$ we require 
$\sum_{u\sim v}m_{uv}=d_v$ instead of $\sum_{u\sim v}m_{uv}=n_v$.  The connection
$\Phi$ is defined as above using $n_v$, however. 

We let $\Omega_{\n,d}$ be the resulting set of multiwebs.
Note that the handshake lemma in this case implies that
$$\sum_{b\in B_{int}}n_b + \sum_{b\in \bB}d_b = \sum_{w\in W_{int}}n_w + \sum_{w\in\bW}d_w$$
in order for $\Omega_{\n,d}$ to be nonempty.

\subsection{Multiweb traces}\label{tracesection}

Material in this section is a generalization of material from \cite{DKS}.
For constant $\n\equiv n$, $n$-webs are used to graphically represent $\SL_n$-invariants in tensor products of representations \cite{sikora_01}.
For general $\n$, $\n$-multiweb traces have no such (global) invariance, but there is a local 
$\SL_{n_v}$-invariance at each vertex $v$; see below.

With boundary vertices of multiplicities $(n_1,d_1),\dots,(n_k,d_k)$, and a multiweb $m \in \Omega_{\n,d}$, 
we will define the trace $\Tr_\Phi(m)$ as an element of
$$(V_1^{\otimes d_1}\otimes\dots\otimes V_k^{\otimes d_k})^*$$
where $V_i=\R^{n_i}$ if the $i$th boundary vertex $v_i$ is black and $V_i=(\R^{n_i})^*$ if $v_i$ is white. 
For a graph without boundary, the trace $\Tr_\Phi(m)$ is a number. The notation $\Tr_\Phi$ emphasizes the
dependence on the connection $\Phi$, but we sometimes just write $\Tr(m)$ to simplify the notation.

To define the multiweb trace we require $\G$ to have the additional structure of 
a total ordering of the edges out of each vertex. Planarity of $\G$ gives a circular ordering to these edges
(by convention, counterclockwise at black vertices and clockwise at white vertices)
which we can augment to a total ordering by specifying a starting half-edge at each vertex. 
Graphically we put a mark, called a \emph{cilium},
in the wedge between the starting and ending edge of the ordering. These cilia then determine the total order at
each vertex. The choice of cilia only affects the sign of the multiweb trace,
and the choice of cilium at a vertex with $n_v$ odd does not affect the multiweb trace at all, so we can and will only define cilia at
vertices of even $n_v$. 

Let $\G$ be a planar bipartite graph with choice of cilia and a connection $\Phi$ as above.
Let $m \in \Omega_{\n,d}$ be an $(\n,d)$-multiweb. To define the multiweb trace $\Tr_\Phi(m)$ it is convenient to first assume that
$m$ is \emph{simple}, that is, all edge multiplicities $m_e$ are $1$ or $0$; we'll discuss higher multiplicity edges below.

Fix an internal black vertex $b$, let $n=n_b$ and $V = \R^{n}$. 
Let $e_1,\dots,e_{n}$ be its adjacent half-edges in order, with multiplicities $m_{e_i}=1$ (we
ignore edges of multiplicity $0$). We place a copy $V^i$ of $V$ on each half-edge $e_i$. 
There is a distinguished element $v_b\in V^1\otimes \dots \otimes V^n$, called the \emph{codeterminant}. In a fixed basis $u_1,\dots,u_{n}$
of $V$,
it is defined as 
\be\label{codet}v_b = \sum_{\sigma\in S_n}(-1)^\sigma u^1_{\sigma(1)}\otimes\dots\otimes u^n_{\sigma(n)}\ee
where $u^i_j$ is the $j$th basis element of $V^i$. 
Note that $v_b$ is independent of $\SL_n$-change of basis of $V$. Also note that if $n$ is even, changing the location of the cilium
by one ``notch" changes the sign of $v_b$, since the signature of the cyclic permutation of length $n$ is $-1$. 

Similarly, for a white interior vertex $w$, we take a copy of $V^* = (\R^{n_w})^*$, the dual space to $\R^{n_w}$, on each half-edge incident to $w$,
and define a codeterminant $v_w\in (V^*)^1\otimes\dots\otimes (V^*)^{n_w}$ in a similar manner. 

The trace $\Tr_\Phi(m)$ is a multilinear
function whose inputs are vectors and covectors at the boundary vertices: at a black boundary vertex $b$ the input is
a vector in $v_b\in (\R^{n_b})^{\otimes d_b}$  and at a white boundary vertex $w$ the input is a vector $v_w\in ((\R^{n_w})^*)^{\otimes d_w}$.
For each vertex of $\G$, we thus have a distinguished element of the corresponding vector space at that vertex:
for internal vertices it is the appropriate codeterminant, and for boundary vertices it is the input vector or covector. 
We take the tensor product of
all these elements to get an element 
$$x \in \bigotimes_{b \in B_{\mathrm{int}}} (\R^{n_b})^{\otimes n_b} \otimes \bigotimes_{w \in W_{\mathrm{int}}} ((\R^{n_w})^*)^{\otimes n_w}\otimes\bigotimes_{b \in \bB} (\R^{n_b})^{\otimes d_b}\otimes\bigotimes_{w \in \bW} ((\R^{n_w})^*)^{\otimes d_w}.$$
Next, we use the connection matrices $\phi_{bw}$, and apply the linear map $\pi = \bigotimes_{bw \in E} \phi_{bw}$ to the tensor factors corresponding to the black vertices,
resulting in $\pi(x) \in \bigotimes_{w \in W} (V^{\otimes n_w} \otimes (V^*)^{\otimes n_w})$. Lastly, we take the trace of this tensor
by contracting at all indices.

This defines the trace for simple multiwebs. Now suppose $m$ has edges of higher multiplicity. Replace each edge $e=bw$ of multiplicity $m_e>1$ with $m_e$ parallel
edges (connecting the same endpoints), extending the linear order at each vertex in the natural way, so that the new
edges are consecutive in the new order. On each new edge use the same mapping $\phi_{bw}$.
This defines a simple multiweb $\hat m$ on this new graph; define 
\be\label{tracemult}\Tr_\Phi(m) = \frac{\Tr_\Phi(\hat m)}{\prod_e m_e!}.\ee

\begin{figure}
\begin{center}\includegraphics[width=1in]{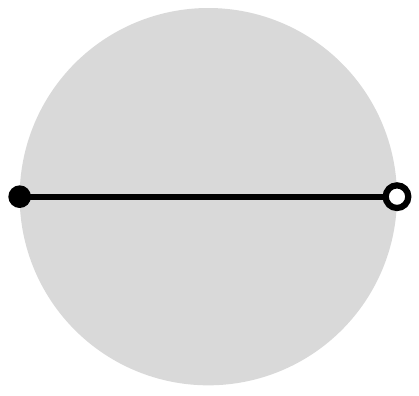} 
\end{center}
\caption{\label{fig:simple_diagram} A simple web.}
\end{figure}

\begin{ex}
The simplest example is depicted in Figure \ref{fig:simple_diagram}.
In this example, the graph is simply one edge connecting two boundary vertices $b,w$ with multiplicities $n_b,n_w$ and $d_b=d_w=1$.
The unique multiweb $m$ has $m_e=1$.
Its trace is $v_w\phi_{bw}v_b$. 
More generally, suppose $d_b=d_w=d$; the unique web has multiplicity $m_e=d$. 
Let $v_b\in(\R^{n_b})^{\otimes d}$ and $v_w\in((\R^{n_w})^*)^{\otimes d}$.
Then the trace is $v_w((\phi_{bw})^{\otimes d})v_b$. In other words, $\Tr_\Phi(m) \in ((\Bbb{R}^{n_b})^{\otimes d} \otimes ((\Bbb{R}^{n_w})^*)^{\otimes d})^*$
is the tensor whose coefficients with respect to the standard basis are the matrix entries of $\phi_{bw}^{\otimes d}$.
\end{ex}

\begin {ex} \label{ex:trivalent_vertex}
Consider the graph $G$ with a single internal 3-valent white vertex and 3 black boundary vertices, shown in Figure \ref{ddex}.
\begin{figure}
\begin{center}\includegraphics[width=3in]{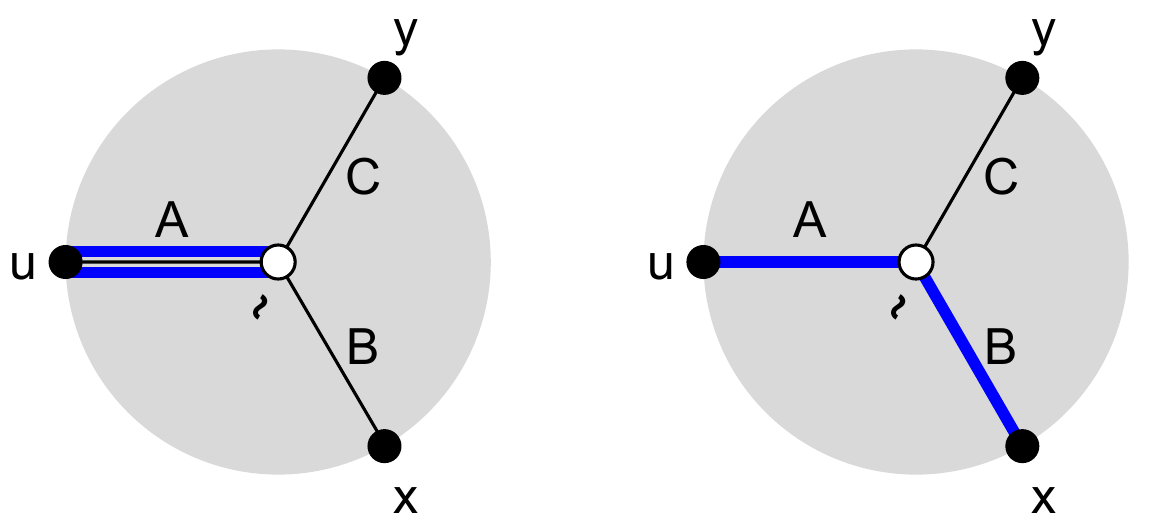} 
\end{center}
\caption{\label{ddex} Two multiwebs.}
\end{figure}
In this example, we take $n_b=n_w= 2$ at all vertices, and the $2 \times 2$ matrices on the edges are $A,B,C$. The input vectors on the boundary are $u,x,y\in(\R^2)^{\otimes d}$,
where $d$ is the corresponding boundary multiplicity. 
The left part of the figure shows an element $m \in \Omega_{2,d}(\mathcal{G})$ where $d_u = 2$, and $d_x=d_y=0$.
The right part of the figure shows an example where $d_u=d_x=1$ and $d_y=0$.
All other $(2,d)$-multiwebs (for all possible choices of $d$) are obtained from these examples by rotation. 
The blue edges indicate the multiplicities. For the left web $u\in\R^2\otimes\R^2$, and web trace is 
$c_a\det A,$ where $c_a$ is the coefficient of the projection of $u$ onto $\R^2\wedge\R^2$.
(That is, letting $\sigma$ denote the exchange of tensor factors of $\R^2\otimes\R^2$, we have
$u_s:=\frac12(u+\sigma u)$ and $u_a:=\frac12(u-\sigma(u))$, and $u_a$ is a multiple $c_a$ of
$e_1\wedge e_2:=e_1\otimes e_2-e_2\otimes e_1$.)
We can see this by splitting the edge into two edges of multiplicity $1$ and applying
(\ref{tracemult}).
Specifically, letting $e_1,e_2$ be the standard basis of $\R^2$, and $f_1,f_2$ the dual basis,
if $u=e_1\otimes e_2$, then $u_a=\frac12e_1\wedge e_2$ and
the trace is 
$$\frac1{2!}\langle v_w|A|u\rangle = \frac12\langle f_1\otimes f_2-f_2\otimes f_1|Ae_1\otimes Ae_2\rangle=\frac12(A_{11}A_{22}-A_{12}A_{21})=\frac12\det A.$$

For the web on the right $u,x\in\R^2$ and the trace is 
$$\langle v_w|A\otimes B|u\otimes x\rangle = \langle f_1\otimes f_2-f_2\otimes f_1|Au\otimes Bx\rangle=Au\wedge Bx,$$
that is, the determinant of the matrix
whose columns are $Au$ and $Bx$. Note that $Au$ comes before $Bx$ due to the location of the cilium, and the fact that the edges are ordered clockwise at a white vertex.
\end{ex}

\subsection{Traces and colorings}\label{trcolorsection}

Let $m\in \Omega_{\n}$. The trace $\Tr(m)$ has an interpretation in terms of edge colorings, as follows. At a vertex $v$ of multiplicity $n_v$ we take a set $\C_v=\{1,2,\dots,n_v\}$ of colors. 
Given a multiweb $m$ and edge $e=bw$ we pick a set $S_e\subset\C_b$ and $T_e\subset \C_w$, both of size $m_e$. 
Around each internal black vertex $b$ we require that the sets $S_{bw}$ for neighbors $w\sim b$ partition $\C_b$,
and likewise at internal white vertices the sets $T_{bw}$ partition $\C_w$. 
We then associate to $b$ a sign $c_b$ depending on the sets 
$\{S_{bw}\}_{w\sim b}$, as follows. Starting from the cilium at $b$, 
order the sets $S_{bw}$ in counterclockwise order, and within each $S_{bw}$ order the colors in natural order.
This gives a permutation of $\C_b$ and $c_b$ is its signature.
Likewise at a white vertex $w$, $c_w$ is the signature of the permutation of the colors around $w$
starting from the cilium in \emph{clockwise} order, with each $T_{bw}$ ordered in natural order.

Let $v$ be a boundary vertex of multiplicity $(n_v,d_v)$, with neighbors $u_1,\dots,u_\ell$. 
Assume for the moment that the cilium at $v$ is in the external face.
Let $u_v\in(\R^{n_v})^{\otimes d_v}$ be a basis vector. It 
corresponds to a map $U_v:[d_v]\to\C_{n_v}$, that is, an ordered $d_v$-tuple of colors at $v$, not necessarily distinct.
Given a multiweb $m$, with multiplicities $m_i$ on edges $vu_i$, we assign to the first edge, $vu_1$, the first $m_1$
colors $U_v(\{1,\dots,m_1\})$, assign to the second edge the next $m_2$ colors $U_v(\{m_1+1,\dots,m_1+m_2\})$, 
and so on. If the cilium at $v$ is not in the external face, we adjust these colors cyclically by the appropriate amount.

The following proposition was proved in \cite{DKS}, Prop. 3.2., for constant $\n$, and no boundary. 
Their proof generalizes directly to the current setting. 

\begin{prop}
For a fixed choice $C_\partial=\{U_v\}_{v\in V_\partial}$ of color sets at the boundary vertices,
the trace of $m$ evaluated at the corresponding basis vectors $\{u_v\}$ is
\be\label{tracecolor} \Tr(m)(\{u_v\}) = \sum_C \prod_v c_v \prod_{e=bw}\det[(\phi_{bw})_{S_e}^{T_e}]\ee
where the sum is over colorings $C$ of the half-edges of $m$ having boundary colors $C_\partial$, assigned as above, 
and $(\phi_{bw})_{S_e}^{T_e}$ is the submatrix
of $\phi_{bw}$ with rows $T_e$ and columns $S_e$.
\end{prop}

As a special case of this proposition, suppose $m$ is simple, that is, $m_e=0$ or $1$
for each edge. Then 
\be\label{tracecolorsimple} \Tr(m)(\{u_v\}) = \sum_C \prod_v c_v \prod_{e=bw}(\phi_{bw})_{t,s}\ee
where the sum is over colorings $C$ of the half-edges of $m$ having boundary colors $C_\partial$ in the correct order, 
and the second product is over edges of $m$, with edge $bw$ getting colors
$S_e=\{s\}$ and $T_e=\{t\}$.
This form in fact implies (\ref{tracecolor}), using the definition (\ref{tracemult}) of trace for edges with $m_e>1$.

\subsection{Kasteleyn matrix}

Given a planar bipartite graph $\G$ with no boundary, with multiplicities $\n$ and connection $\Phi$, 
let $L$ be a choice of cilia at vertices of even multiplicity. 
We choose a $\pm1$-connection $\eps=\{\eps_{e}\}_{e\in E}$
on $\G$, the \emph{Kasteleyn connection}, or  \emph{Kasteleyn signs}, with monodromy $(-1)^{\ell+1+k}$ around each face,
where $2\ell$ is the length of the face (its number of edges) and $k$ is the number of cilia in that face. 
Recall that we do not specify cilia at vertices with $n_v$ odd, so $k$ only counts cilia at vertices with $n_v$ even.

We define a matrix $K=K(\Phi)$, the \emph{Kasteleyn matrix}, as follows. 
If $w\sim b$, $K_{wb}=\eps_{wb}\phi_{bw}$ and otherwise $K_{wb}=0$.
Note that this generalizes the usual Kasteleyn rule of \cite{Kenyon.lectures}
(which corresponds to the
case when all faces have an even number of cilia; in this case the sign only depends on the face length mod $4$; see also \cite{DKS}). 

The following lemma will be used below.
\begin{lemma}\label{loopsign}
For any loop $\gamma$ in $\G$, the monodromy of $\eps$ around $\gamma$ is
$(-1)^{L/2+1+k+n_{\mathrm{int}}}$, where $L$ is the length of $\gamma$,
$k$ is the number of cilia of vertices along $\gamma$ which point into $\gamma$,
and $n_{\mathrm{int}}$ is the sum of multiplicities of vertices
strictly enclosed by $\gamma$. 
\end{lemma}

\begin{proof}
This is true by definition when $\gamma$ is a single face and follows by induction on the number of faces enclosed by $\gamma$: suppose we subdivide a face $f$ into two faces by adding in a path between two vertices $u,v$ of $f$,
so that $f$ becomes two faces $f_1,f_2$ (see Figure \ref{subfaces}). 
\begin{figure}
\begin{center}
\includegraphics[width=2.6in]{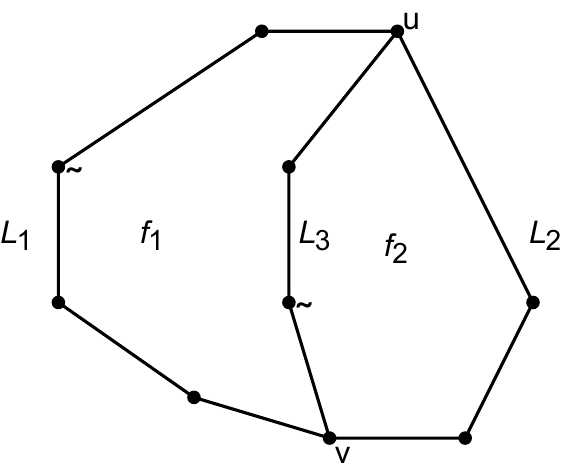}    \end{center}
\caption{\label{subfaces}Subdividing a face $f$ into two faces $f_1,f_2$.}
\end{figure} 
Let $L_1$ be the length (number of edges)
of the part $\rho_1$ of $f$ in common with $f_1$, and let $L_2$ be the length of the part $\rho_2$ of $f$ in common with $f_2$. Let $L_3$ be the length of the new path.
Let $k_1$ be the number of cilia along $\rho_1$ (pointing into $f$), and $k_2$ that for $\rho_2$.
Let $k_3$ be the cilia along the new path pointing into $f_1$, and $k_4$ the cilia along the new path pointing
into $f_2$.
We have
$$(-1)^{\frac{L_1+L_3}2+1+k_1+k_3}(-1)^{\frac{L_2+L_3}2+1+k_2+k_4}=(-1)^{\frac{L_1+L_2}2+1+k_1+k_2+(L_3+k_3+k_4+1)}.$$
It remains to identify $L_3+k_3+k_4+1\bmod 2$ with $n_{int}\bmod 2$.
Each of the $L_3-1$ vertices along the new path either has a cilium (in $f_1$ or $f_2$, and is thus of even multiplicity) 
or is of odd multiplicity. The sum $k_3+k_4$ is the number of even-multiplicity vertices along the new path, so 
$L_3-k_3-k_4-1=n_{int}\bmod 2$.
\end{proof}

We then let $\tilde K =\tilde K(\Phi)$ be obtained from $K$ be replacing each entry $K_{wb}$,
with the corresponding $n_w\times n_b$ block of scalars (in particular $0$ entries of $K$ become
blocks of $0$s).

\begin {ex} \label{ex:Kexample}
Consider the graph $\mathcal{G}$ pictured in Figure \ref{fig:Kexample} where the $\phi_{bw}$ are $A,B,\dots,G$. 
Note that $\eps_{bw}=+1$ everywhere except on the edge $33$ (hence the sign in the label $-F$).
We have
\begin {equation} \label{eqn:Kasteleyn_example}
    K = \begin{pmatrix}
        A & D & 0 \\
        B & C & G \\
        E & 0 & -F
    \end{pmatrix}, ~~~~
    \widetilde{K} = \begin{pmatrix}
        a   & d_1    & d_2    & 0      & 0 \\
        b_1 & c_{11} & c_{12} & g_{11} & g_{12} \\
        b_2 & c_{21} & c_{22} & g_{21} & g_{22} \\
        b_3 & c_{31} & c_{32} & g_{31} & g_{32} \\
        e   & 0      & 0      & -f_1   & -f_2
    \end{pmatrix}.
\end {equation}
\begin{figure}
\begin{center}
\includegraphics[width=2.6in]{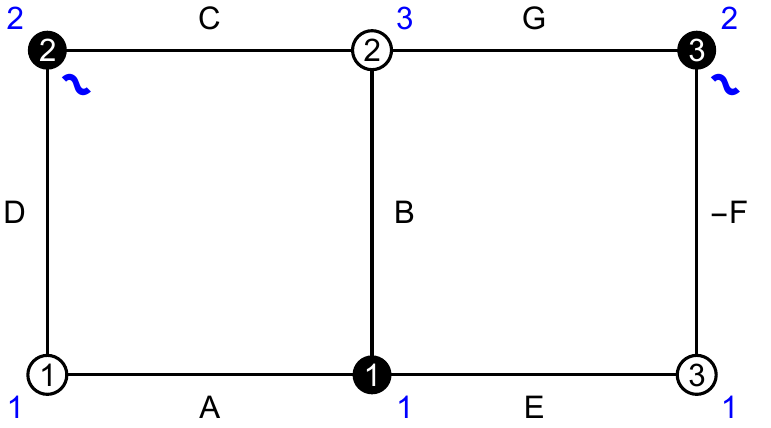}    \end{center}
\caption{\label{fig:Kexample}Example of $\G$ with multiplicities, signs and cilia. The labels on the vertices are the indices, and blue numbers are the multiplicities.}
\end{figure}
\end {ex}

\subsection{Trace theorem}

When $\n\equiv 2$, a $2$-multiweb, also known as a double-dimer cover, is a collection of loops and doubled edges. 
If the connection takes values in $\SL_2$, the web trace is the product of the matrix traces
of the monodromies of the connection around these loops (with either orientation). In \cite{kenyon_14}, it was shown that the web trace of a general bipartite graph with an $\SL_2$-connection
is given by a generalization of the Kasteleyn determinant formula. This was generalized in \cite{DKS} both to $n > 2$ and to the case of connection matrices in $\mathrm{GL}_n$ and even $M_n(\R)$.
We further generalize this result (with a similar proof) to the case of varying $n$ as follows.

\begin {thm}\label{main}
Let $\G$ be a bipartite planar graph without boundary and with choice of cilia $L$. 
Let $\n:V\to\Z_{>0}$, and
choose a connection $\Phi$. 
Let $K$ be an associated Kasteleyn matrix. Then
\[ \det(\widetilde{K}) = \pm\sum_{m\in\Omega_{\n}} \Tr_\Phi(m).  \]
\end {thm}

\begin{proof}
Let $\G_{\n}$ be obtained from $\G$ by replacing each vertex $v$ with $n_v$ vertices $v_1,\dots,v_{n_v}$, and each edge $bw$ by a complete bipartite graph $K_{n_b,n_w}$. 
On edge $b_jw_i$ lying over edge $bw$ of $\G$, put (scalar) weight $\eps_{bw}(\phi_{bw})_{i,j}$, 
where $\eps_{bw}$ is the Kasteleyn sign from $bw$ in $\G$ and $(\phi_{bw})_{i,j}$ is the 
$i,j$-entry of $\phi_{bw}$. 
Note that $\det\tilde K$ is a (signed, weighted) sum of dimer covers of $\G_{\n}$.

There is a natural projection map $\G_{\n}\to\G$. Each dimer cover of $\G_{\n}$ projects to an $\n$-multiweb of $\G$,
and we can group the terms in $\det\tilde K$ according to their corresponding multiwebs:
$$\det\tilde K = \sum_{m\in\Omega_\n}\sum_{\sigma\in m}(-1)^\sigma \tilde K_\sigma.$$
Here $\sigma\in m$ means the dimer cover $\sigma$ projects to $m$, and $\tilde K_\sigma$ is the
corresponding product of entries of $\tilde K$. 

It remains to show that, up to a global sign $s$,
\be\label{traceterm1}s \Tr(m) = \sum_{\sigma\in m}(-1)^\sigma \tilde K_\sigma= 
\sum_{\sigma\in m}(-1)^\sigma\prod_{\tilde w}\tilde K_{\tilde w,\sigma(\tilde w)}.\ee

We assume first that $m$ is simple, that is, all edge multiplicities $m_e$ are $0$ or $1$. If $m$ is not simple we temporarily modify the graph by replacing each edge on which $m_e>1$ with $m_e$ parallel edges; see below.

Let $w\in\G$ and $\tilde w\in\G_{\n}$ lying over $w$. Then $\sigma(\tilde w)=\tilde b$ is a black vertex of $\G_{\n}$ lying over some neighbor $b$ of $w$ in $\G$. Under the identification
of vertices over $v$ with colors at $v$, this $\sigma$ thus colors the edge $e=wb$, with one
color from $\C_w$ and one color from $\C_b$. 
These colors on edge $wb$ also correspond to a matrix entry $\phi_{ij}$. 
Thus (\ref{traceterm1}) agrees with (\ref{tracecolorsimple}), once we verify the agreement of signs.  We need to check that 
\be\label{signs}(-1)^\sigma\prod_{bw\in m}\eps_{bw}\prod_v c_v=s,\ee
a constant sign.  

Suppose we exchange colors at two circularly adjacent half-edges at $v$ (not separated by the cilium at $v$).
Then both $(-1)^\sigma$ and $c_v$ change signs. Thus all colorings of a given simple 
multiweb $m$
have the same left-hand side of (\ref{signs}). 

It remains to check the sign change when
we change multiwebs.
But first let us check what happens when we glue parallel edges together to make non-simple multiwebs. 
If simple multiweb $m$ has $k$ parallel edges
between $w$ and $b$, consider all colorings of (the half edges of) these edges with colors from size-$k$ subsets $S_e\subset \C_w$ and $T_e\subset\C_b$. For fixed $S_e,T_e$
there are $(k!)^2$ such colorings;
summing over all $\sigma$ involving such colorings
gives a contribution of $k!\det ((\phi_{wb})_{S_e}^{T_e})$. When we glue these parallel edges of $m$
into a single edge of multiplicity $k$, and divide by $k!$, the contribution from this 
edge from these colors is $\det ((\phi_{wb})_{S_e}^{T_e})$, which agrees with the corresponding term in (\ref{tracecolor}).
So at this point we can assume $m$ is a general (nonsimple) multiweb.

Now we need to check the sign change when
we change multiwebs. 
Let $\gamma$ be a simple closed path in $\G$ with vertices $b_1,w_1,\dots,b_\ell,w_\ell$ in counterclockwise order. 
Suppose multiweb $m\in\Omega_{\n}$ 
has positive multiplicities on edges $b_iw_i$ of $\gamma$. Let $m'$
be obtained from $m$ by decreasing the multiplicities on edges $b_iw_i$ by $1$
and increasing the multiplicities by $1$ on the other edges $b_iw_{i-1}$ (indices cyclic). 
We show below that these operations preserve (\ref{signs}) and connect the set of all $\n$-multiwebs; this will complete the 
proof.

By changing colors if necessary (we already dealt with color swaps above) we can assume
that both half edges of $b_iw_i$ contain color $1$, the first color at each of $b,w$. 
We then move these 
colors $1$ to the edges
$b_iw_{i-1}$, increasing their multiplicity. 
Note that $m'\in\Omega_{\n}$ as well. Let us compare the sign (\ref{signs}) of $m$ with $m'$.
Suppose edge $b_iw_i$ has multiplicity $m_i$, and let $q(b_i)$ be the total multiplicity
of edges $b_iw$ of $m$ with $w$ strictly enclosed by $\gamma$. 
If the cilium at $b_i$ is also in the region enclosed by $\gamma$,
then when moving color $1$ the sign change of $c_{b_i}$ is $(-1)^{m_{i}+q(b_i)}$:
color $1$ is the first color of edge $b_iw_i$, so moves past $m_i-1+q(b_i)$ other colors and also past the cilium, to become the first color of $b_iw_{i-1}$.
Similarly if the cilium at $b$ is not enclosed by $\gamma$ the sign change of $c_{b_i}$ is $(-1)^{m_i+q(b_i)-1}$.
Likewise at $w_i$ the sign change is $(-1)^{m_i+q(w_i)}$ if the cilium at $w_i$ is enclosed,
and otherwise is $(-1)^{m_i+q(w_i)-1}$. 
Therefore the total sign change in the $c_v$, taking the product around all vertices of $\gamma$,
is $(-1)^{k+Q}$,
where $k$ is the number of cilia at the vertices of $\gamma$ which are enclosed by $\gamma$,
and $Q=\sum_{v\in\gamma}q(v)$. Note that, by the handshake lemma, 
$(-1)^Q=(-1)^{n_{\mathrm{int}}}$ where $n_{\mathrm{int}}$ 
is the total multiplicity of  vertices strictly enclosed by $\gamma$. 

This sign change $(-1)^{k+n_{\mathrm{int}}}$ is exactly compensated by the monodromy of the Kasteleyn signs $\eps$ around the face,
so that (\ref{signs}) remains unchanged. More precisely, 
the change in the sign is due to the change in the signature of $\sigma$, which is $(-1)^{\ell+1}$ (since $\sigma$ changes by multiplication by an $\ell$-cycle) 
times the monodromy around the face of the Kasteleyn signs,
which is $(-1)^{\ell+1+k+n_{\mathrm{int}}}$ by Lemma \ref{loopsign}, times the contribution from the cilia of $(-1)^k$. 
The result is no net sign change. 

It remains to show that $\Omega_{\n}$ is connected by such ``loop moves". 
We can interpret each $m\in\Omega_{\n}$ as a flow $[m]\in\Z^E$, with value $m_e$ on edge
$e$ directed from black to white. Then given $m,m'\in\Omega_{\n}$, the difference 
$[m]-[m']$ is a divergence-free flow. It can be decomposed into a collection of non-crossing 
oriented cycles, that is, disjoint cycles on the associated ribbon surface. We can
change $m$ to $m'$ on each cycle using the operation described above.
\end{proof}

\begin {ex}
    Consider again the case from Example \ref{ex:Kexample}, pictured in Figure \ref{fig:Kexample}. It is easy to see by inspection
    that there are only 3 elements of $\Omega_{\n}(\mathcal{G})$, pictured below:
    \begin {center}
    \begin {tikzpicture}
        \draw (0,0) grid (2,1);
        \draw[red, line width=1.5] (0,0) -- (1,0);
        \draw[red, line width=1.5] (0,1) -- (1,1);
        \draw[red, line width=1.5] (0,1.1) -- (1,1.1);
        \draw[red, line width=1.5] (1,1) -- (2,1);
        \draw[red, line width=1.5] (2,0) -- (2,1);

        \foreach \x/\y in {0/1, 1/0, 2/1}
            \draw[fill=black] (\x,\y) circle (0.06);
        \foreach \x/\y in {0/0, 1/1, 2/0}
            \draw[fill=white] (\x,\y) circle (0.06);

        \draw (1,-0.5) node {$m_1$};

        \begin {scope}[shift={(4,0)}]
            \draw (0,0) grid (2,1);
            \draw[red, line width=1.5] (0,0) -- (0,1);
            \draw[red, line width=1.5] (0,1) -- (1,1);
            \draw[red, line width=1.5] (1,1.1) -- (2,1.1);
            \draw[red, line width=1.5] (1,1) -- (2,1);
            \draw[red, line width=1.5] (2,0) -- (1,0);

            \foreach \x/\y in {0/1, 1/0, 2/1}
                \draw[fill=black] (\x,\y) circle (0.06);
            \foreach \x/\y in {0/0, 1/1, 2/0}
                \draw[fill=white] (\x,\y) circle (0.06);

            \draw (1,-0.5) node {$m_2$};
        \end {scope}

        \begin {scope}[shift={(8,0)}]
            \draw (0,0) grid (2,1);
            \draw[red, line width=1.5] (0,0) -- (0,1);
            \draw[red, line width=1.5] (0,1) -- (2,1);
            \draw[red, line width=1.5] (1,1) -- (1,0);
            \draw[red, line width=1.5] (2,0) -- (2,1);

            \foreach \x/\y in {0/1, 1/0, 2/1}
                \draw[fill=black] (\x,\y) circle (0.06);
            \foreach \x/\y in {0/0, 1/1, 2/0}
                \draw[fill=white] (\x,\y) circle (0.06);

            \draw (1,-0.5) node {$m_3$};
        \end {scope}
    \end {tikzpicture}
    \end {center}
    In this case, one can see that $\sum_{i=1}^3 \Tr(m_i)$ realizes the expansion of $\det(\widetilde{K})$ along the first column.
    Indeed, if $\widetilde{K}^{ij}$ is the matrix obtained by removing row $i$ and column $j$, then $\Tr(m_1) = a \det(\widetilde{K}^{11})$,
    $\Tr(m_2) = e \det(\widetilde{K}^{51})$, and $\Tr(m_3) = b_1\det(\widetilde{K}^{21}) + b_2\det(\widetilde{K}^{31}) + b_3\det(\widetilde{K}^{41})$.
\end {ex}

\section {Grassmannians and Boundary Measurements} \label{sec:grassmannians}

In this section, we review Postnikov's parameterization of the positive Grassmannian using networks, which he calls
\emph{plabic graphs} (the word plabic is an abbreviation of ``planar bi-colored'') \cite{postnikov_06}. These are planar graphs, drawn in a disk,
with some vertices on the boundary. Vertices are drawn in two colors (black and white).

The \emph{Grassmannian} $\Gr_{k,n}$ is the set of $k$-dimensional subspaces of an $n$-dimensional vector space.
It can be viewed concretely as the quotient of the set of $k \times n$ matrices, of rank $k$, modulo row operations (i.e. left multiplication by $\mathrm{GL}_k$).
The $k$-dimensional subspace represented by a matrix is simply its row span. The \emph{Pl\"{u}cker coordinates} are the minors $\Delta_I$,
where $I \in \binom{[n]}{k}$ is a $k$-element subset of $[n] = \{1,2,\dots,n\}$, and $\Delta_I$ is the determinant of the $k \times k$ submatrix using 
the columns from $I$. The \emph{non-negative Grassmannian} $\Grnn_{k,n}$ is the set of points whose Pl\"{u}cker coordinates are
all non-negative. 

Postnikov showed that $\Grnn_{k,n}$ has a cell decomposition whose cells, called \emph{positroid cells}, are indexed by
equivalence classes of plabic graphs. For each ``reduced" plabic graph\footnote{A plabic graph is \emph{reduced} if it has the minimal number of faces required to parameterize the appropriate cell; see \cite{postnikov_06, gk_13}.}, 
the assignment of weights (positive real numbers) to edges (or faces) of the graph 
gives a parameterization of the corresponding positroid cell. This parameterization was called the \emph{boundary measurement map} in \cite{postnikov_06},
and we will briefly describe it now.

An orientation of a plabic graph is called a \emph{perfect orientation} if every black vertex has in-degree 1, and every white vertex has out-degree 1.
The following construction uses a choice of perfect orientation, but the end result will not depend on the choice. The chosen orientation will make
some of the boundary vertices sources, and others sinks. Suppose $k$ is the number of sinks (and $n-k$ the number of sources), and $I \subset [n]$ the set of sinks. 
The corresponding $k \times n$ matrix $X$, called the \emph{boundary measurement matrix}, is defined as follows. The rows will correspond to the sinks in
the set $I$, and the columns will correspond to all boundary vertices. The submatrix in columns $I$ will simply be the identity matrix.
For all other entries, in row $i$ and column corresponding to boundary vertex $j$, the entry of the matrix will be the generating function for directed
paths from source $j$ to sink $i$:
\[ X_{ij} = (-1)^{N_{ij}} \sum_{p \colon j \to i} (-1)^{w(p)} \mathrm{wt}(p) \]
Here, the sum is over all directed paths from $j$ to $i$, $N_{ij}$ is the number of sinks on the boundary between $i$ and $j$, 
$w(p)$ is the winding number of the path, and $\mathrm{wt}(p)$ is simply the product of the edge weights along the path.
If the graph is acyclic, then the factors of $(-1)^{w(p)}$ are irrelevant. When there are cycles, these signs ensure that the rational expressions
obtained by summing the geometric series are subtraction-free.

\begin {ex} \label{ex:boundary_measurement}
    Figure \ref{fig:boundary_measurements} shows a network, and two choices of perfect orientations. As in the previous section, we use the
    convention that the preferred orientation is from black-to-white, and all our edge weights are with respect to this convention.
    Therefore when a perfect orientation assigns some edges a white-to-black direction, we invert the weight on that edge. 
    The boundary measurement matrices of the two orientations are 
    \[ 
        X_1 = \begin{pmatrix} 
            1 & 0 & -\frac{c}{d} & -\frac{1}{d} \\
            0 & 1 & b + \frac{ac}{d} & \frac{a}{d}
        \end{pmatrix} \quad \text{ and } \quad
        X_2 = \begin{pmatrix} 
            1 & \frac{c}{ac+bd} & 0 & -\frac{b}{ac+bd} \\
            0 & \frac{d}{ac+bd} & 1 & \frac{a}{ac+bd}
        \end{pmatrix} 
    \]
    Note that the entries of $X_1$ are Laurent polynomials, while the entries of $X_2$ are genuine rational expressions. This comes from the directed
    cycle in the second perfect orientation.
    Note however that the two matrices are equivalent up to left action by $\mathrm{GL}_2$,
    which can be seen by the fact that they have the same $2 \times 2$ minors (up to a common factor of $\frac{d}{ac+bd}$).

    \begin {figure}[h!]
    \centering
    \begin {tikzpicture}
        \draw (-0.5,0) -- (2.5,0) -- (2.5,2) -- (-0.5,2) -- cycle;

        \draw (0.5,0.5) -- (1.5,0.5) -- (1.5,1.5) -- (0.5,1.5) -- cycle;
        \draw (-0.5,0.5) -- (0.5,0.5);
        \draw (-0.5,1.5) -- (0.5,1.5);
        \draw (2.5,0.5) -- (1.5,0.5);
        \draw (2.5,1.5) -- (1.5,1.5);
        \draw[fill=white] (0.5,0.5) circle (0.06);
        \draw[fill=white] (1.5,1.5) circle (0.06);
        \draw[fill=black] (0.5,1.5) circle (0.06);
        \draw[fill=black] (1.5,0.5) circle (0.06);
        \draw[fill=white] (0,1.5) circle (0.06);
        \draw[fill=white] (2,0.5) circle (0.06);

        \draw (-0.5,1.5) node[left]  {$1$};
        \draw (-0.5,0.5) node[left]  {$2$};
        \draw (2.5,0.5)  node[right] {$3$};
        \draw (2.5,1.5)  node[right] {$4$};

        \draw[fill=black] (-0.5,1.5) circle (0.06);
        \draw[fill=black] (-0.5,0.5) circle (0.06);
        \draw[fill=black] (2.5,1.5)  circle (0.06);
        \draw[fill=black] (2.5,0.5)  circle (0.06);

        \draw (0.5,1)  node[left]  {$a$};
        \draw (1,0.55) node[below] {$b$};
        \draw (1.5,1)  node[right] {$c$};
        \draw (1,1.45) node[above] {$d$};

        \begin {scope}[shift={(4.4,0)}]
            \draw (-0.5,0) -- (2.5,0) -- (2.5,2) -- (-0.5,2) -- cycle;

            \draw[-latex] (1.5,1.5)  -- (0.56,1.5);
            \draw[-latex] (0.5,1.5)  -- (0.5,0.56);
            \draw[-latex] (1.5,0.5)  -- (1.5,1.44);
            \draw[-latex] (1.5,0.5)  -- (0.56,0.5);
            \draw[-latex] (0.5,0.5)  -- (-0.44,0.5);
            \draw[-latex] (0.5,1.5)  -- (0.06,1.5);
            \draw[-latex] (0,1.5)    -- (-0.44,1.5);
            \draw[-latex] (2.5,0.5)  -- (2.06,0.5);
            \draw[-latex] (2,0.5)    -- (1.56,0.5);
            \draw[-latex] (2.5,1.5)  -- (1.56,1.5);

            \draw[fill=white] (0.5,0.5) circle (0.06);
            \draw[fill=white] (1.5,1.5) circle (0.06);
            \draw[fill=black] (0.5,1.5) circle (0.06);
            \draw[fill=black] (1.5,0.5) circle (0.06);
            \draw[fill=white] (0,1.5) circle (0.06);
            \draw[fill=white] (2,0.5) circle (0.06);

            \draw[fill=black] (-0.5,1.5) circle (0.06);
            \draw[fill=black] (-0.5,0.5) circle (0.06);
            \draw[fill=black] (2.5,1.5)  circle (0.06);
            \draw[fill=black] (2.5,0.5)  circle (0.06);

            \draw (-0.5,1.5) node[left]  {$1$};
            \draw (-0.5,0.5) node[left]  {$2$};
            \draw (2.5,0.5) node[right] {$3$};
            \draw (2.5,1.5) node[right] {$4$};

            \draw (0.5,1)  node[left]  {$a$};
            \draw (1,0.55) node[below] {$b$};
            \draw (1.5,1)  node[right] {$c$};
            \draw (1,1.45) node[above] {$d^{-1}$};
        \end {scope}

        \begin {scope}[shift={(8.8,0)}]
            \draw (-0.5,0) -- (2.5,0) -- (2.5,2) -- (-0.5,2) -- cycle;

            \draw[-latex] (0.5,1.5)  -- (1.44,1.5);
            \draw[-latex] (0.5,0.5)  -- (0.5,1.44);
            \draw[-latex] (1.5,0.5)  -- (0.56,0.5);
            \draw[-latex] (1.5,1.5)  -- (1.5,0.56);
            \draw[-latex] (-0.5,0.5) -- (0.44,0.5);
            \draw[-latex] (0.5,1.5)  -- (0.06,1.5);
            \draw[-latex] (0,1.5)    -- (-0.44,1.5);
            \draw[-latex] (1.5,0.5)  -- (1.94,0.5);
            \draw[-latex] (2,0.5)    -- (2.44,0.5);
            \draw[-latex] (2.5,1.5)  -- (1.56,1.5);

            \draw[fill=white] (0.5,0.5) circle (0.06);
            \draw[fill=white] (1.5,1.5) circle (0.06);
            \draw[fill=black] (0.5,1.5) circle (0.06);
            \draw[fill=black] (1.5,0.5) circle (0.06);
            \draw[fill=white] (0,1.5) circle (0.06);
            \draw[fill=white] (2,0.5) circle (0.06);

            \draw[fill=black] (-0.5,1.5) circle (0.06);
            \draw[fill=black] (-0.5,0.5) circle (0.06);
            \draw[fill=black] (2.5,1.5)  circle (0.06);
            \draw[fill=black] (2.5,0.5)  circle (0.06);

            \draw (-0.5,1.5) node[left]  {$1$};
            \draw (-0.5,0.5) node[left]  {$2$};
            \draw (2.5,0.5)  node[right] {$3$};
            \draw (2.5,1.5)  node[right] {$4$};

            \draw (0.5,1)  node[left]  {$\frac{1}{a}$};
            \draw (1,0.55) node[below] {$b$};
            \draw (1.5,1)  node[right] {$\frac{1}{c}$};
            \draw (1,1.45) node[above] {$d$};
        \end {scope}
    \end {tikzpicture}
    \caption {(left) A plabic graph for $\Gr_{2,4}$, (middle) an acyclic perfect orientation, (right) a cyclic perfect orientation.}
    \label {fig:boundary_measurements}
    \end {figure}
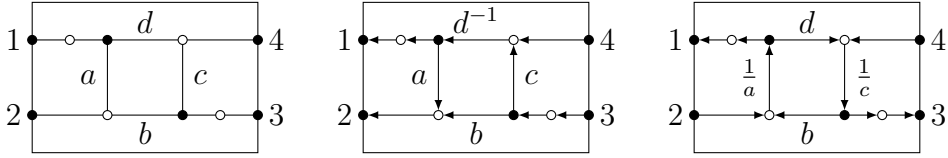
\end {ex}

In \cite{psw_09}, an interpretation of the boundary measurements in terms of dimer covers was given. Suppose the graph is bipartite,
and all boundary vertices are black (this can be arranged by some local transformations which do not change the boundary measurements,
specifically by inserting some 2-valent vertices).
An \emph{almost perfect matching} is a subset of edges where every internal vertex is incident to exactly one edge (but boundary vertices
are allowed to be unmatched). It can be shown that for a given graph, there is a number $k$ (called the \emph{excedence}) such that
every almost perfect matching uses exactly $k$ boundary vertices. This is the same $k$ as above, since there is a simple bijection
between almost perfect matchings and perfect orientations: the edges of the matching are the edges directed white-to-black in the perfect orientation.
Then the sinks of the perfect orientation are exactly the boundary vertices used in this matching.

If $M$ is an almost perfect matching, let $\partial(M) \subset [n]$ be the set of boundary vertices used by $M$. The following result
can be seen in \cite{psw_09, KenyonWilson, lam_15}. 

\begin {thm} \label{thm:pluckers_and_dimers}
    Let $X$ be the boundary measurement matrix corresponding to a chosen perfect orientation. Then
    \[ \Delta_I(X) = \sum_{M \colon \partial(M) = I} \mathrm{wt}(M) \]
    In other words, the Pl\"{u}cker coordinate $\Delta_I$ is the generating function for almost perfect matchings with boundary set $I$.
\end {thm}

The boundary measurement matrix can be given the following
concrete interpretation in terms of the Kasteleyn matrix. Suppose $\G$ is bipartite, all
boundary vertices are black,
and let $K$ be a (typically not square) Kasteleyn matrix for $\G$. Let $m$ be an almost perfect matching, and 
$W'$ the set of white vertices paired in $m$ to a boundary vertex.
Then after renumbering the vertices, 
$K$ has block form $K=\begin{pmatrix}A&B\\C&D\end{pmatrix}$
where $D$ has rows indexed by $W\setminus W'$ and columns indexed by $B_{\mathrm{int}}$,
the internal black vertices. Note that $D$ is a square matrix since its vertices are paired by $m$. Then the boundary measurement matrix, up to a fixed choice of 
gauge, is the Schur reduction
$X=A-BD^{-1}C$; see \cite{KenyonWilson, KenyonShi}. 
The matrix $X$ will be totally nonnegative
after a gauge transformation; a choice of Kasteleyn signs which yields totally nonnegative $X\in\Grnn_{k,n}$ is one where, for each $i=1,\dots,n-1$, 
the product of signs along the boundary 
path between adjacent boundary vertices $i$ and $i+1$ is $(-1)^{\ell/2+1}$ where $\ell$ is the length of the path.

\begin {ex}
    Consider the situation from Example \ref{ex:boundary_measurement}, pictured in Figure \ref{fig:boundary_measurements}.
    Figure \ref{fig:plucker_dimers} illustrates Theorem \ref{thm:pluckers_and_dimers} in this example, showing the various almost perfect matchings.
    The Pl\"{u}cker coordinates are 
    \[ \Delta_{12} = d, \quad \quad \Delta_{13} = ac+bd, \quad \quad \Delta_{14} = a \]
    \[ \Delta_{23} = c, \quad \quad \Delta_{24} = 1, \quad \quad \Delta_{34} = b \]
    These agree with the minors of $X_1$ up to a common factor of $d$, and they agree with the minors of $X_2$ up to a common factor of $\frac{1}{ac+bd}$.

    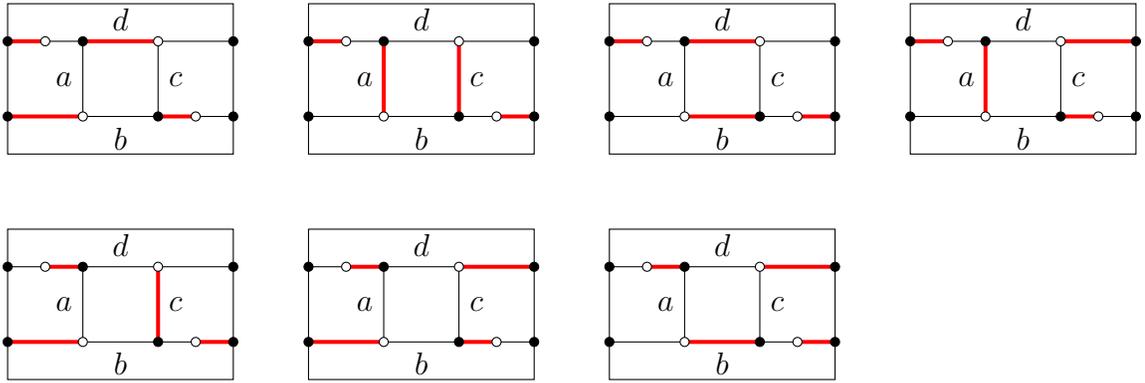
\begin {figure}[h!]
    \centering
    \begin {tikzpicture}


        \foreach \x/\y in {0/0, 4/0, 8/0, 12/0, 0/-3, 4/-3, 8/-3} {
            \draw (\x,\y) -- (\x+3,\y) -- (\x+3,\y+2) -- (\x,\y+2) -- cycle;

            \draw (\x+0,\y+0.5) -- (\x+3,\y+0.5);
            \draw (\x+0,\y+1.5) -- (\x+3,\y+1.5);
            \draw (\x+1,\y+0.5) -- (\x+1,\y+1.5);
            \draw (\x+2,\y+0.5) -- (\x+2,\y+1.5);
        }



        \draw[red, line width = 1.5] (0,1.5) -- (0.5,1.5);
        \draw[red, line width = 1.5] (0,0.5) -- (1,0.5);
        \draw[red, line width = 1.5] (1,1.5) -- (2,1.5);
        \draw[red, line width = 1.5] (2,0.5) -- (2.5,0.5);


        \begin {scope}[shift={(4,0)}]
            \draw[red, line width = 1.5] (0,1.5) -- (0.5,1.5);
            \draw[red, line width = 1.5] (1,1.5) -- (1,0.5);
            \draw[red, line width = 1.5] (2,1.5) -- (2,0.5);
            \draw[red, line width = 1.5] (2.5,0.5) -- (3,0.5);
        \end {scope}

        \begin {scope}[shift={(8,0)}]
            \draw[red, line width = 1.5] (0,1.5) -- (0.5,1.5);
            \draw[red, line width = 1.5] (1,1.5) -- (2,1.5);
            \draw[red, line width = 1.5] (1,0.5) -- (2,0.5);
            \draw[red, line width = 1.5] (2.5,0.5) -- (3,0.5);
        \end {scope}


        \begin {scope}[shift={(12,0)}]
            \draw[red, line width = 1.5] (0,1.5) -- (0.5,1.5);
            \draw[red, line width = 1.5] (2,1.5) -- (3,1.5);
            \draw[red, line width = 1.5] (1,0.5) -- (1,1.5);
            \draw[red, line width = 1.5] (2,0.5) -- (2.5,0.5);
        \end {scope}


        \begin {scope}[shift={(0,-3)}]
            \draw[red, line width = 1.5] (0,0.5) -- (1,0.5);
            \draw[red, line width = 1.5] (2.5,0.5) -- (3,0.5);
            \draw[red, line width = 1.5] (0.5,1.5) -- (1,1.5);
            \draw[red, line width = 1.5] (2,0.5) -- (2,1.5);
        \end {scope}


        \begin {scope}[shift={(4,-3)}]
            \draw[red, line width = 1.5] (0,0.5) -- (1,0.5);
            \draw[red, line width = 1.5] (2,1.5) -- (3,1.5);
            \draw[red, line width = 1.5] (0.5,1.5) -- (1,1.5);
            \draw[red, line width = 1.5] (2,0.5) -- (2.5,0.5);
        \end {scope}


        \begin {scope}[shift={(8,-3)}]
            \draw[red, line width = 1.5] (2,1.5) -- (3,1.5);
            \draw[red, line width = 1.5] (2.5,0.5) -- (3,0.5);
            \draw[red, line width = 1.5] (1,0.5) -- (2,0.5);
            \draw[red, line width = 1.5] (0.5,1.5) -- (1,1.5);
        \end {scope}


        \foreach \x/\y in {0/0, 4/0, 8/0, 12/0, 0/-3, 4/-3, 8/-3} {
            \draw[fill=white] (\x+1,\y+0.5)   circle (0.06);
            \draw[fill=white] (\x+2,\y+1.5)   circle (0.06);
            \draw[fill=black] (\x+1,\y+1.5)   circle (0.06);
            \draw[fill=black] (\x+2,\y+0.5)   circle (0.06);
            \draw[fill=white] (\x+0.5,\y+1.5) circle (0.06);
            \draw[fill=white] (\x+2.5,\y+0.5) circle (0.06);

            \draw[fill=black] (\x, \y+1.5)   circle (0.06);
            \draw[fill=black] (\x, \y+0.5)   circle (0.06);
            \draw[fill=black] (\x+3, \y+0.5) circle (0.06);
            \draw[fill=black] (\x+3, \y+1.5) circle (0.06);

            \draw (\x+1,  \y+1)   node[left]  {$a$};
            \draw (\x+1.5,\y+0.5) node[below] {$b$};
            \draw (\x+2,  \y+1)   node[right] {$c$};
            \draw (\x+1.5,\y+1.5) node[above] {$d$};
        }
    \end {tikzpicture}
    \caption {The almost perfect matchings of the graph from Figure \ref{fig:boundary_measurements}}
    \label {fig:plucker_dimers}
    \end {figure}
\end {ex}

Lastly, we will describe an inverse to the boundary measurement map, at least in a special case. Among the positroid cells,
there is one of maximal dimension, which a dense open set in $\Grnn_{k,n}$. It is called the \emph{totally positive}
Grassmannian, and denoted $\Grtp_{k,n}$. The graphs used to parameterize this top-dimensional
cell correspond to the standard cluster algebra structure on the coordinate ring of the Grassmannian. They can be obtained as follows.
Start by drawing a $k \times (n-k)$ rectangular array of vertices. The vertices along the left and bottom edges will be connected to the boundary.
Direct all edges in the rectangular grid up and to the left. This is pictured in the left part of Figure \ref{fig:top_cell_graph}. 
To turn this into a perfectly oriented plabic graph, we can split all 4-valent vertices into two 3-valent vertices. The resulting graph
will be perfectly oriented, and we can color the vertices black and white accordingly. Also, to make the graph more uniform (and bipartite), we can insert some 2-valent vertices.
This is pictured in the right part of Figure \ref{fig:top_cell_graph}.

\begin {figure}[h!]
\centering
\begin {tikzpicture}
    \draw (0,0) -- (5,0) -- (5,4) -- (0,4) -- cycle;

    \foreach \x in {1,2,3,4} {
        \foreach \y in {1,2,3} {
            \draw[fill=black] (\x,\y) circle (0.06);
        }
    }

    \foreach \x in {0,1,2,3} {
        \foreach \y in {1,2,3} {
            \draw[latex-] (\x+0.06,\y) -- (\x+1,\y);
        }
    }

    \foreach \x in {1,2,3,4} {
        \foreach \y in {1,2,3} {
            \draw[latex-] (\x,\y-0.06) -- (\x,\y-1);
        }
    }

    \foreach \x in {1,2,3,4} {
        \draw[fill=black] (\x,0) circle (0.06);
    }

    \foreach \y in {1,2,3} {
        \draw[fill=black] (0,\y) circle (0.06);
    }


    \begin {scope}[shift={(7,0)}]
        \draw (0,0) -- (5,0) -- (5,4) -- (0,4) -- cycle;

        \foreach \x/\y in {1/2.732, 2/2.732, 3/2.732, 1.5/1.866, 2.5/1.866, 3.5/1.866, 2/1, 3/1} {
            \draw[latex-] (\x+0.03,\y+0.03) -- ($(\x,\y) + (0.45,0.239)$);
        }

        \foreach \x/\y in {1.5/3.021, 2.5/3.021, 3.5/3.021, 1/2.155, 2/2.155, 3/2.155, 4/2.155, 1.5/1.289, 2.5/1.289, 3.5/1.289} {
            \draw[latex-] ($(\x,\y) + (0.03,-0.03)$) -- ($(\x,\y) + (0.45,-0.239)$);
        }

        \foreach \x/\y in {1/2.732, 2/2.732, 3/2.732, 4/2.732, 1.5/1.866, 2.5/1.866, 3.5/1.866, 4.5/1.866} {
            \draw[latex-] (\x,\y-0.05) -- ($(\x,\y) - (0,0.527)$);
        }

        \draw[latex-] ($(4.5,1.289) - (0.45,0.239)$) -- ($(4.5,1.289) + (-0.03,-0.03)$);

        \draw[latex-] (0,2.732) -- (0.95, 2.732);
        \draw[latex-] (0.05,2.155) -- (0.5, 2.155);
        \draw[latex-] (0.55,2.155) -- (1, 2.155);
        \draw[latex-] (0.05,1.289) -- (0.5, 1.289);
        \draw[latex-] (0.55,1.289) -- (1.5, 1.289);
        \draw[latex-] (2,1) -- (2, 0);
        \draw[latex-] (3,1) -- (3, 0);
        \draw[latex-] (4,1) -- (4, 0);

        \draw[latex-] (4.5,0.55) -- (4.5, 0);
        \draw[latex-] (4.5,1.239) -- (4.5, 0.6);

        \draw[fill=white] (0.5,2.155) circle (0.06);
        \draw[fill=white] (0.5,1.289) circle (0.06);
        \draw[fill=white] (4.5,0.6)   circle (0.06);

        \foreach \x/\y in {1.5/3.021, 2.5/3.021, 3.5/3.021, 1/2.155, 2/2.155, 3/2.155, 4/2.155, 1.5/1.289, 2.5/1.289, 3.5/1.289, 4.5/1.289} {
            \draw[fill=black] (\x,\y) circle (0.05);
        }

        \foreach \x/\y in {1/2.732, 2/2.732, 3/2.732, 4/2.732, 1.5/1.866, 2.5/1.866, 3.5/1.866, 4.5/1.866, 2/1, 3/1, 4/1} {
            \draw[fill=white] (\x,\y) circle (0.05);
        }

        \foreach \x in {2,3,4,4.5} {
            \draw[fill=black] (\x,0) circle (0.05);
        }

        \foreach \y in {1.289,2.155,2.732} {
            \draw[fill=black] (0,\y) circle (0.05);
        }
    \end {scope}
\end {tikzpicture}
\caption {The plabic graph representing the top-dimensional cell $\Grtp_{3,7}$ in $\Grnn_{3,7}.$}
\label {fig:top_cell_graph}
\end {figure}

\section{Multiweb probability models}

A \emph{multiweb probability model}
is a certain probability measure on $\Omega_{\n}$, the set of multiwebs on a given bipartite planar graph
(without boundary), with a connection having the property that the trace is a positive function
on $\Omega_{\n}$. In this case we define $Pr(m) = \frac{\Tr(m)}Z$
where $Z=\sum_{m\in\Omega_{\n}}\Tr(m)$ is the normalizing constant, called the \emph{partition function}.  By Theorem \ref{main}, $Z$ is the determinant of an associated Kasteleyn matrix.

We give some examples here. Some results of \cite{KOS} are used; we refer the reader there or to \cite{Kenyon.lectures}
for the relevant background.

\subsection{Free-fermionic six-vertex model}\label{6vsection}

\begin{figure}
\begin{center}
\includegraphics[width=2in]{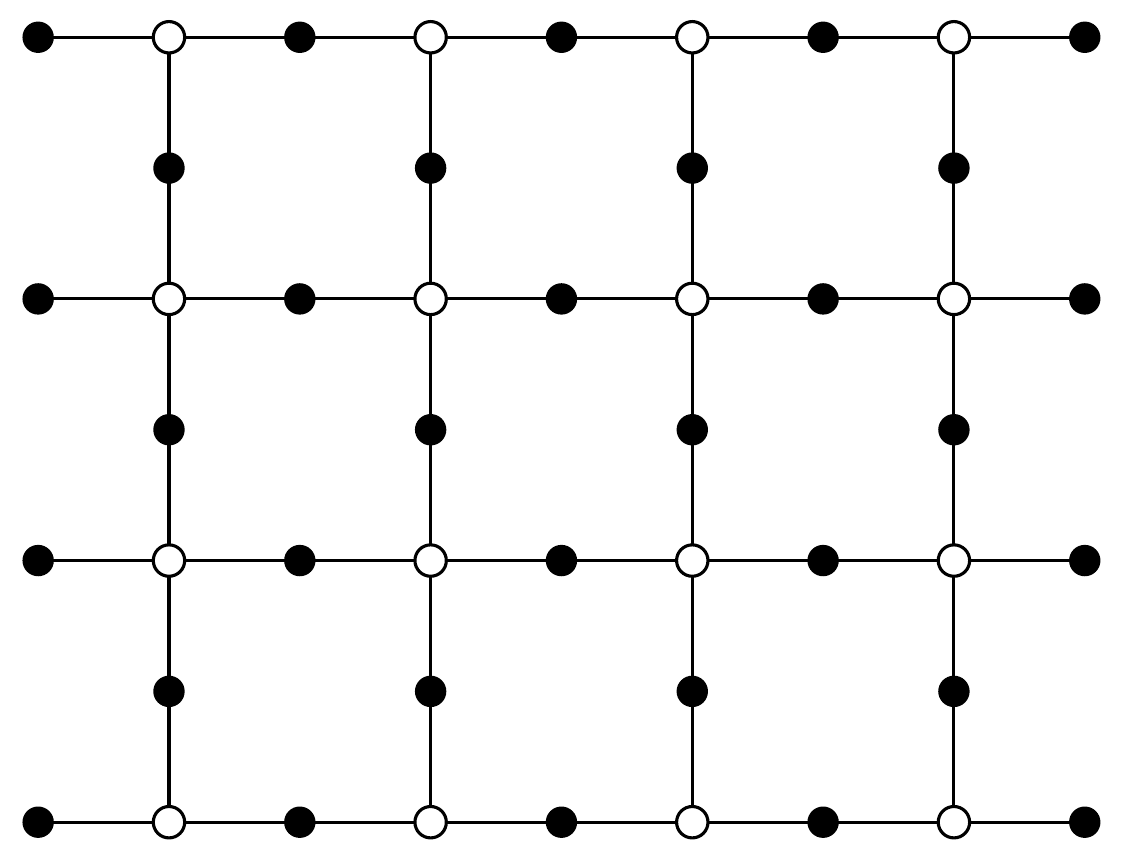}\hskip1cm\includegraphics[width=1.5in]{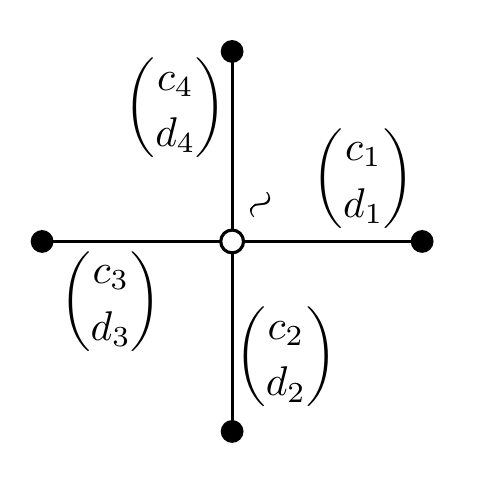}
\end{center}
\caption{\label{6Vmatrices}Graph for the $6V$ multiweb model, and connection matrices around a white vertex.}
\end{figure}

Consider the $6$-vertex model discussed in the introduction.
We work with the square ice version.
Let $\G$ be the graph obtained from $\Z^2$ by putting a vertex in the center of each edge.
Color these new vertices black and the original vertices white (Figure \ref{6Vmatrices}, left).
The set of square ice configurations is in bijection with the space $\Omega_{\n}$ of ${\bf n}$-multiwebs 
on $\G$, where $n_w=2$ at white vertices and $n_b=1$ at black vertices.

We put $2\times 1$ matrices $\phi_{bw}$ on the edges, and cilia northeast of each white vertex as shown in Figure \ref{6Vmatrices}, right (we don't need to specify cilia at black vertices). 
Then for any configuration $m$ its trace is the product, over all white vertices $w$,
of the ``molecule" $m_w$ there. The trace of the molecule $m_w$ is a $2\times 2$ minor
of the $2\times4$ Grassmannian matrix 
$$M_w=\begin{pmatrix} c_1&c_2&c_3&c_4\\d_1&d_2&d_3&d_4\end{pmatrix}.$$ 
Specifically the trace of $m_w$ is the minor $\Delta_{i,j}$ where $i,j$
are the edges of the molecule attached to $w$, indexed clockwise starting from the cilium. If we choose each $M_w$ to be in $\Grtp_{2,4}$,
then all traces $\Tr(m_w)$ are positive and $\Tr(m)=\prod_{w\in W}\Tr(m_w)$ is positive.
These traces define a probability measure on $\Omega_{\n}$, that is, on the space of square ice configurations, with $Pr(m) = \frac{\Tr(m)}{Z}$,
where $Z$, the partition function, is $Z=\det K(\Phi)$. 

The six vertex weights $a_1,\dots,c_2$ defined by $M_w$ are the six $2\times2$ minors of
$M_w$. As such they satisfy a relation,
the Pl\"ucker relation $a_1a_2+b_1b_2-c_1c_2=0$.  This equation defines the so-called
\emph{free-fermionic subvariety} of the six vertex model.  It is well-known that the model
can be solved by determinantal methods in this case, see e.g. \cite{Baxter}.

A particularly symmetric choice is $a_1=a_2=b_1=b_2=1$ and $c_1=c_2=\sqrt{2}$.
In this case the square ice model on the $n\times n$ grid graph with ``domain wall boundary conditions", that is, for the graphs
shown in Figure \ref{6V} or Figure \ref{6Vmatrices}, left,
is equivalent to domino tilings of the Aztec diamond, see \cite{EKLP2}. 

For the model on the whole plane defined by $M_w\in\Grtp_{2,4}$, we can compute the free energy per site $F$
as the Mahler measure of the characteristic polynomial $P(z,w)=\det \tilde K(z,w)$ (see \cite{KOS})
where, taking a fundamental domain consisting of a single white vertex and its east and north black neighbors,
$K(z,w)$ is the $1\times 2$ matrix 
$$K(z,w)=\begin{pmatrix}\begin{pmatrix}c_1\\d_1\end{pmatrix}-z\begin{pmatrix}c_3\\d_3\end{pmatrix}&
\begin{pmatrix}c_2\\d_2\end{pmatrix}+w\begin{pmatrix}c_4\\d_4\end{pmatrix}\end{pmatrix}.$$
Here the minus sign is a Kasteleyn sign. We have 
$$\tilde K(z,w) = \begin{pmatrix}c_1-c_3z&c_2+c_4w\\d_1-d_3z&d_2+d_4w\end{pmatrix}$$
and 
$$P(z,w) =\det\tilde K=  \Delta_{12}+z \Delta_{23}+w \Delta_{14}-wz\Delta_{34}.$$
The free energy per site is
$$F = \frac1{(2\pi)^2}\int_{0}^{2\pi}\int_0^{2\pi}\log P(e^{i\theta},e^{i\phi})\,d\theta\,d\phi.$$

\subsection{20-vertex model}\label{20vsection}

The $20$-vertex model is defined similarly to the $6$-vertex model, but on the triangular lattice.
A configuration consists in an orientation of the edges so that at each vertex the outdegree and indegree are $3$. 
There is an equivalent ``triangular ice" version, for example ammonia ice, $\text{NH}_3$, with a nitrogen atom at each vertex and 
three hydrogen atoms bonded to it along the edges.

\begin{figure}
\begin{center}
\includegraphics[width=2in]{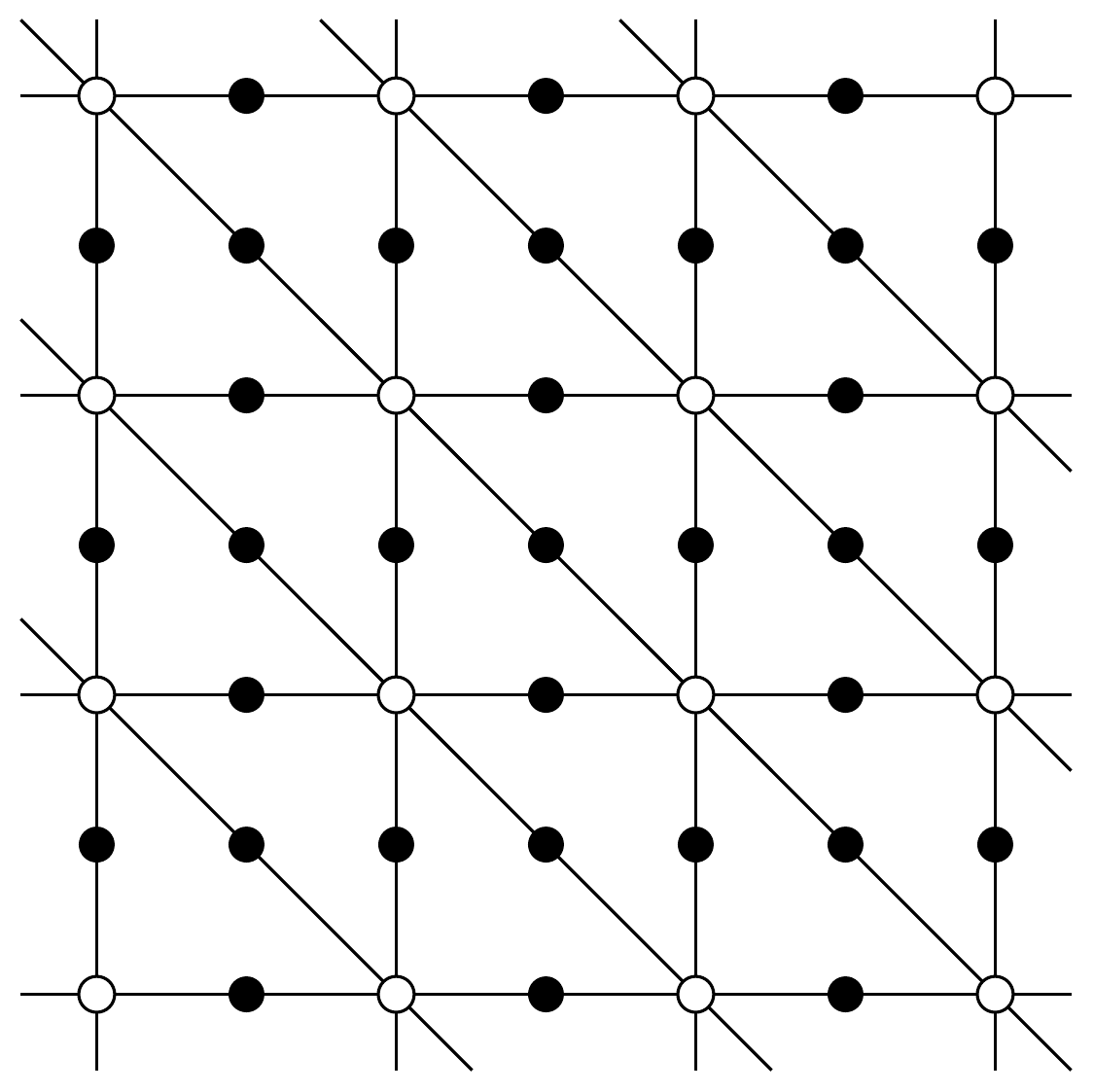}\end{center}
\caption{\label{20Vgraph}}
\end{figure}

To represent the $20$-vertex model as a web model,
let $\G$ be the graph obtained
from the triangular grid by putting a vertex in the center of each edge.
Color these new vertices black and the original vertices white. 
See Figure \ref{20Vgraph}.
The set of ice configurations is in bijection with the space $\Omega_{\n}$ of ${\bf n}$-multiwebs 
on $\G$, where $n_w=3$ at white vertices and $n_b=1$ at black vertices.

We put $3\times 1$ matrices $\phi_{bw}$ on the edges, and periodic cilia, for example east-southeast of each white vertex.
The trace of the molecule $m_w$ is a $3\times 3$ minor
of the $3\times6$ Grassmannian matrix $M_w$: the minor $\Delta_{i,j,k}$ where $i,j,k$
index the edges of the molecule attached to $w$. 
As above we choose each $M_w$ to be in $\Grtp_{3,6}$.
then $\Tr(m)=\prod_{w\in W}\Tr(m_w)$ is positive.
These traces define a probability measure on $\Omega_{\n}$, that is, on the space of 
triangular ice configurations, with $Pr(m) = \frac{\Tr(m)}{Z}$,
where $Z=\det K(\Phi)$. 

The $20$ vertex weights $w_1,\dots,w_{20}$ defined by $M_w$ are the maximal minors of
$M_w$. As such they satisfy certain Pl\"ucker relations.  These equation define a 
$9$-dimensional \emph{free-fermionic subvariety} of the 20-vertex model.

\begin{figure}
\begin{center}
\includegraphics[width=3in]{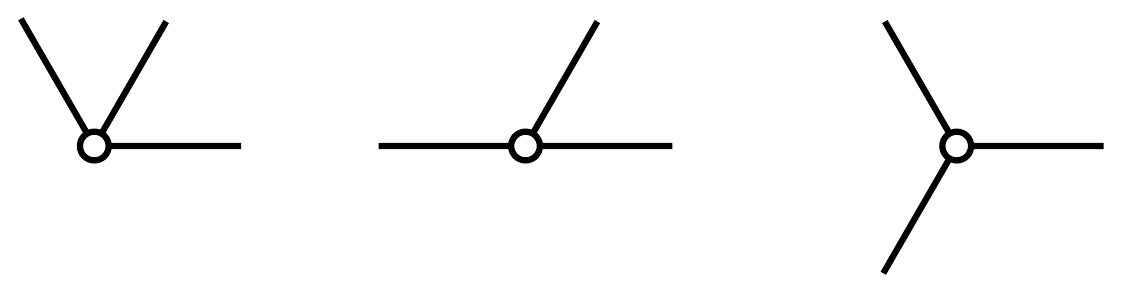}\end{center}
\caption{\label{20Vmolecules}At the symmetric point, the three types of molecules have weights in ratio $1:2:3$.}
\end{figure}

As in the square ice model, there is a circularly symmetric choice 
of weights, arising from the unique circularly symmetric point in $\Grtp_{3,6}$ described
in \cite{Karp}. 
Up to a global scale factor this $6$-fold symmetric version of the model has 
$$M=\begin{pmatrix}1 & 1 & 1 & 1 & 1 & 1 \\
\cos0&  \cos\frac{\pi}3 & \cos\frac{2\pi}3 &\cos\frac{3\pi}3 &  \cos\frac{4\pi}3 & \cos\frac{5\pi}3 \\
\sin0 &  \sin\frac{\pi}3 & \sin\frac{2\pi}3 &\sin\frac{3\pi}3 &  \sin\frac{4\pi}3 & \sin\frac{5\pi}3 \end{pmatrix}
$$
\old{$$M=\begin{pmatrix}1 & \frac{1}{\omega } & \frac{1}{\omega ^2} &
   \frac{1}{\omega ^3} & \frac{1}{\omega ^4} &
   \frac{1}{\omega ^5} \\
 1 & 1 & 1 & 1 & 1 & 1 \\
 1 & \omega  & \omega ^2 & \omega ^3 & \omega ^4 & \omega ^5\end{pmatrix}$$
where $\omega=e^{2\pi i/6}$. (This $M$ is not totally positive but becomes so if any row is multiplied by $i$.)}
If we embed $\G$ symmetrically, so that faces are equilateral,
then each molecule gets a weight proportional to the product, over pairs
of hydrogen atoms, of sines of
the half-angles between vectors connecting the central atom to the arms. As such the 
weights of the $3$ possible types of molecules in Figure \ref{20Vmolecules} are in ratio $1:2:3$. 

We can compute the free energy for the $20$-vertex model on the 
plane by choosing a fundamental
domain consisting of a single white vertex and the three black neighbors east, north, and northwest from it.
Then 
$K(z,w)$ is the $1\times 3$ matrix 
$$K(z,w)=\begin{pmatrix}\begin{pmatrix}c_1\\d_1\\e_1\end{pmatrix}+z\begin{pmatrix}c_4\\d_4\\e_4\end{pmatrix}&
\begin{pmatrix}c_2\\d_2\\e_2\end{pmatrix}+w\begin{pmatrix}c_5\\d_5\\e_5\end{pmatrix}&
\begin{pmatrix}c_3\\d_3\\e_3\end{pmatrix}+\frac{w}{z}\begin{pmatrix}c_6\\d_6\\e_6\end{pmatrix}\end{pmatrix}.$$
We have 
$$\tilde K(z,w) = \begin{pmatrix}c_1+c_4z&c_2+c_5w&c_3+c_6w/z\\d_1+d_4z&d_2+d_5w&d_3+d_6w/z\\
e_1+e_4z&e_2+e_5w&e_3+e_6w/z\end{pmatrix},$$
and 
$$P(z,w) =\det\tilde K=  \Delta_{123}+z \Delta_{234}+zw \Delta_{345}-w(\Delta_{135}+\Delta_{246})+w^2\Delta_{456}+\frac{w^2}{z}\Delta_{156}+\frac{w}{z}\Delta_{126}.$$
The free energy per site is
$$F = \frac1{(2\pi)^2}\int_{0}^{2\pi}\int_0^{2\pi}\log P(e^{i\theta},e^{i\phi})\,d\theta\,d\phi.$$

\begin{figure}
\begin{center}
\includegraphics[width=2in]{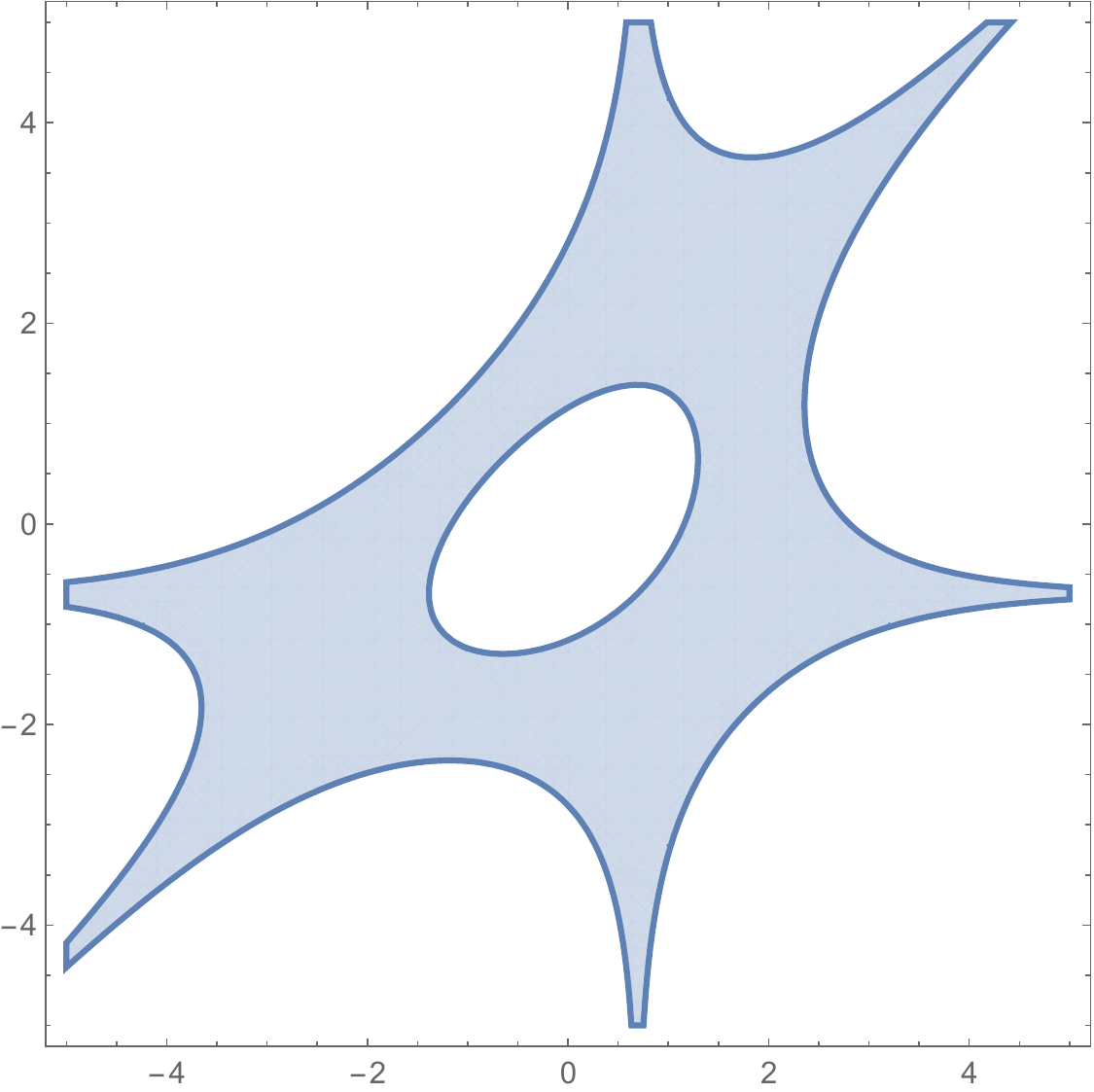}
\end{center}
\caption{\label{amoeba20}The amoeba for $P$ in the $3$-fold symmetric $20$-vertex model.}
\end{figure}

\begin{figure}
\begin{center}
\includegraphics[width=3in]{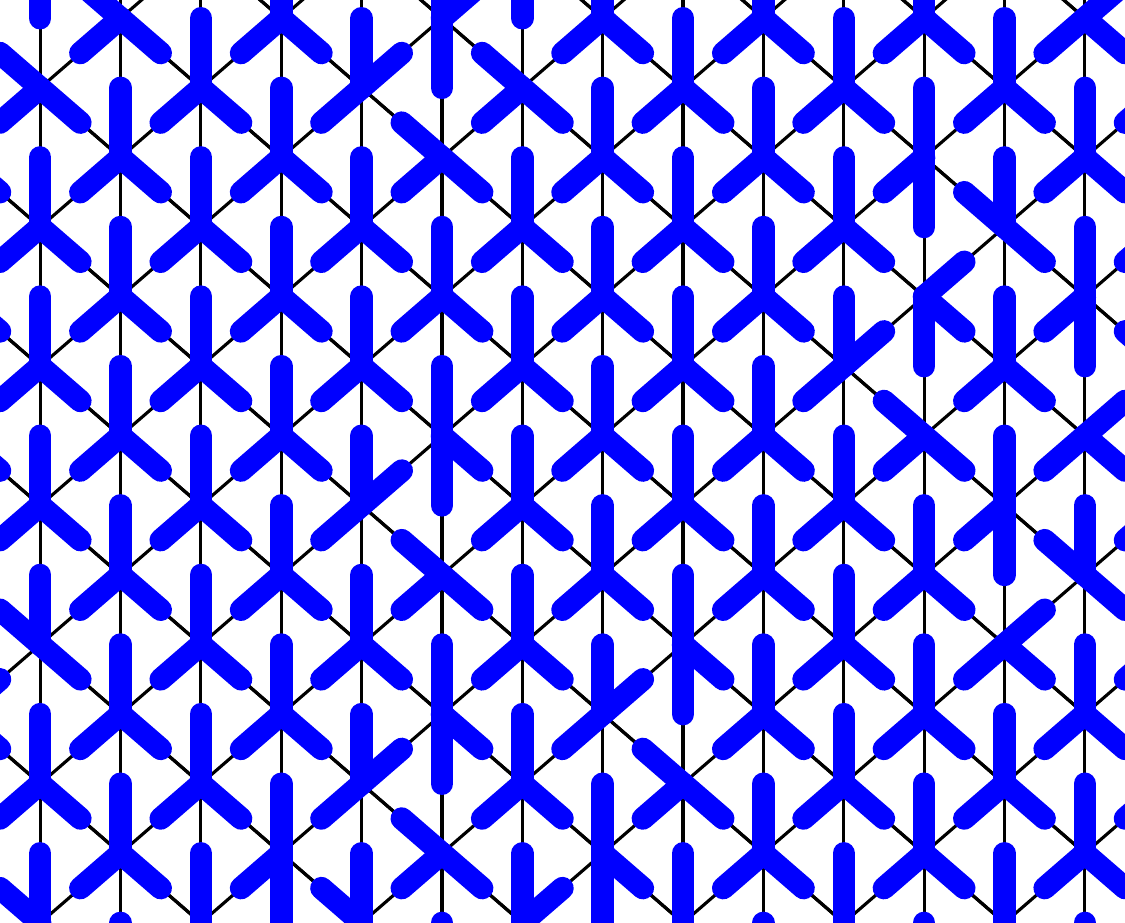}
\end{center}
\caption{\label{locked}Sample configuration from the maximal Gibbs measure for the $20$-vertex model when $a=2$.}
\end{figure}

As a specific example, take
$$M=\begin{pmatrix}a & 1 & a & 1 & a & 1 \\
a\cos0&  \cos\frac{\pi}3 & a\cos\frac{2\pi}3 &\cos\frac{3\pi}3 &  a\cos\frac{4\pi}3 & \cos\frac{5\pi}3 \\
a\sin0 &  \sin\frac{\pi}3 & a\sin\frac{2\pi}3 &\sin\frac{3\pi}3 &  a\sin\frac{4\pi}3 & \sin\frac{5\pi}3 \end{pmatrix}
$$
for a parameter $a>0$. This leads to a model with $3$-fold (but not $6$-fold) rotational symmetry; edges from a white vertex
in directions 
of cube roots of unity have weight $a$ and edges in direction which are cube roots of $-1$ have weight $1$. 
A molecule has a weight $1,2,$ or $3$ as before multiplied by its edge weights. 
We have 
$$P(z,w) = a^2+a z+a^2zw+a w^2+a^2\frac{w^2}z+a\frac{w}{z}-(3+3a^2)w$$
which is easier to understand in ``Newton polygon" form
$$P(z,w) = \begin{matrix}a^2\frac{w^2}z&a w^2\\a\frac{w}{z}&-(3+3a^2)w&a^2 wz\\&a^2&a z\end{matrix}.$$
This model has a ``gas phase" when $a\ne1$, see Figure \ref{amoeba20}: a bounded complementary component of the amoeba of $P$. The maximal entropy Gibbs measure on 
configurations on $\G$ is this gas phase; it has the property that correlations between local arrangements decay exponentially quickly in the distance. A typical configuration has molecules
locked in a periodic pattern except for localized perturbations, as in Figure \ref{locked}.
Interestingly, as $a\to 1$ the gaseous phase shrinks
to a point; at $a=1$ we are in the $6$-fold symmetric case and this is a so-called ``liquid phase'':
the correlations decay only polynomially.

\section {Realization of Higher-Rank Connections by Scalar Networks}\label{scalarsection}

We show here that to a bipartite planar graph $\G$, with the data of multiplicities
$\n$ and connection $\Phi$,
we can associate
a new planar bipartite graph $\hat\G$ with an $\R^*$-connection for which there is a local,
measure-preserving surjective mapping $\Psi:\Omega_1(\hat\G)\to\Omega_{\n}(\G)$.

We begin by illustrating our construction with a worked example. 
We then discuss the general construction below.

\subsection{Worked example} \label{sec:worked_example}

We give an explicit realization of a $\GL_2$ connection on a trivalent bipartite graph $\G$,
by an $\R^*$-connection on a new graph $\hat\G$.

\begin{figure}
\begin{center}
\includegraphics[width=2in]{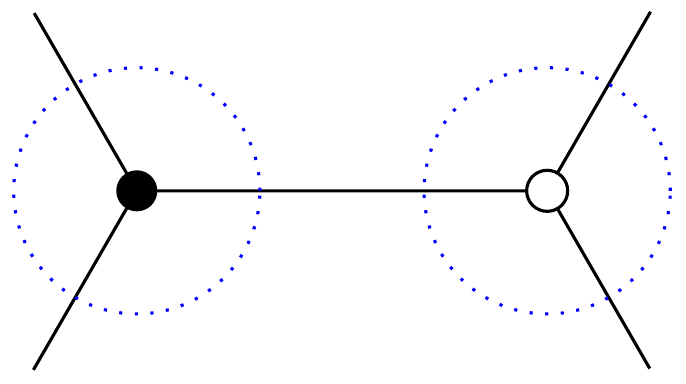}\hskip1cm\includegraphics[width=2in]{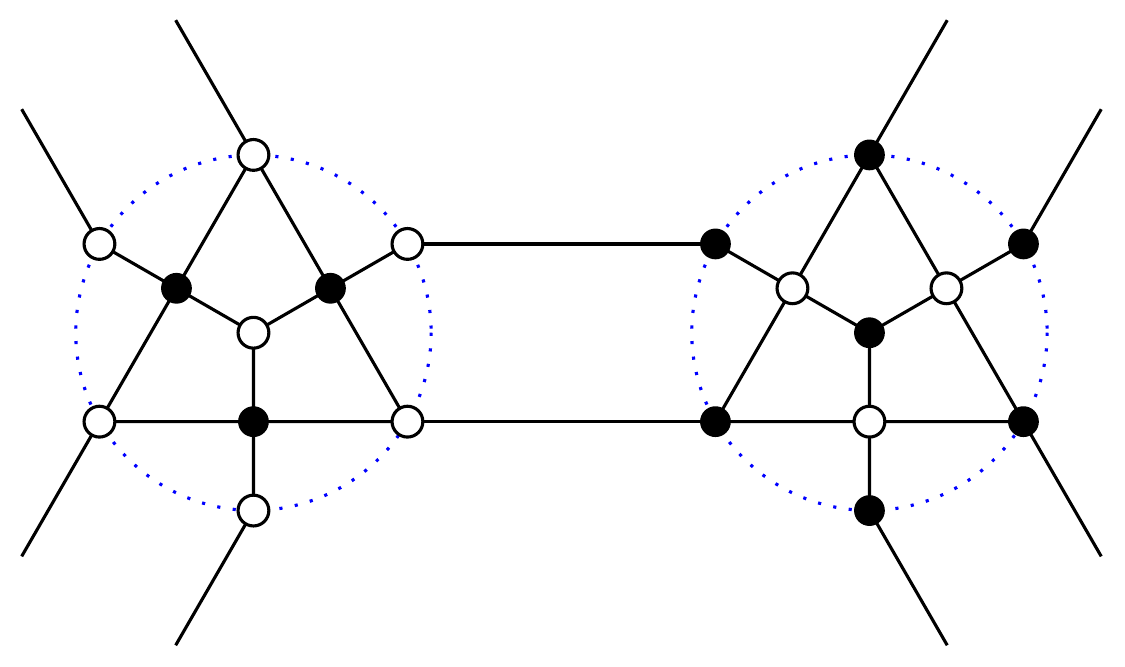}
\end{center}
\caption{\label{hatG2}(Part of) a trivalent graph $\G$ and its ``scalarization" $\hat\G$. }
\end{figure}

\begin{figure}
    \begin{center}
        \includegraphics[width=2in]{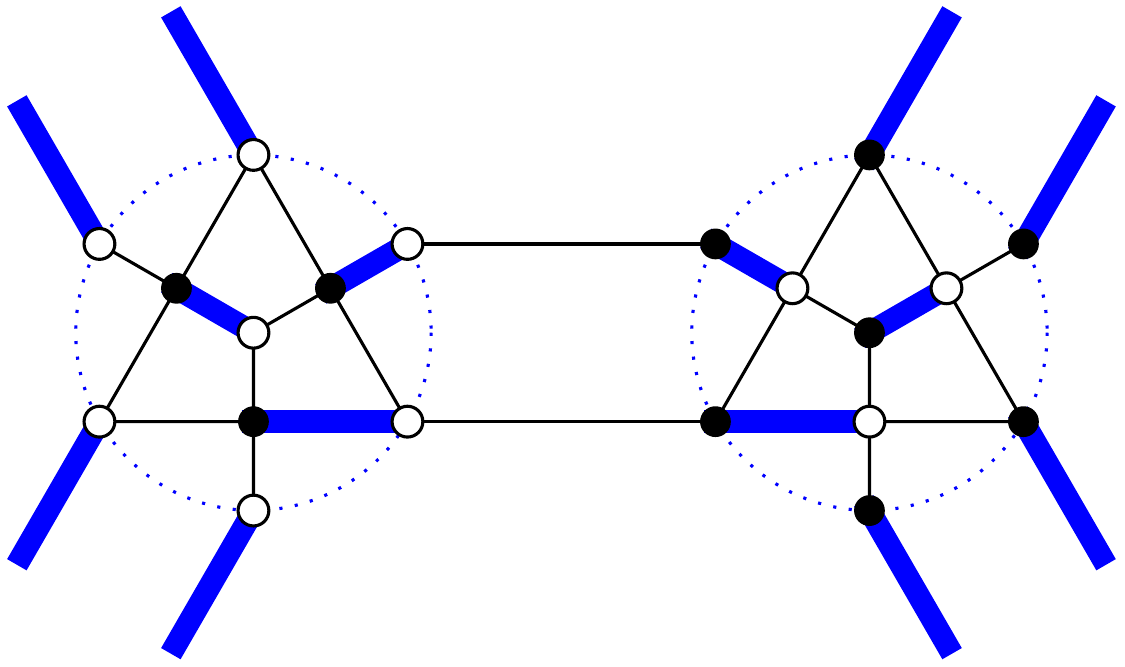}\hskip1cm\includegraphics[width=1.4in]{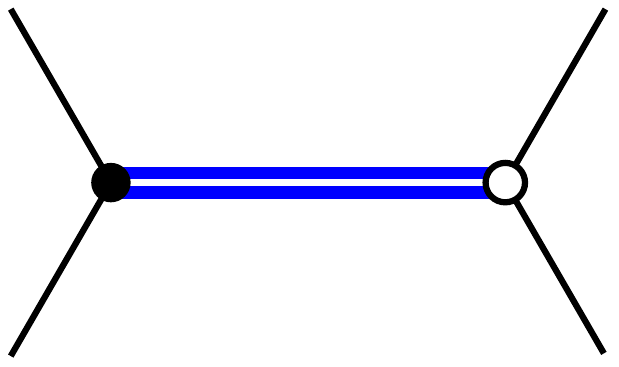}\\
        \vskip.5cm
        \includegraphics[width=2in]{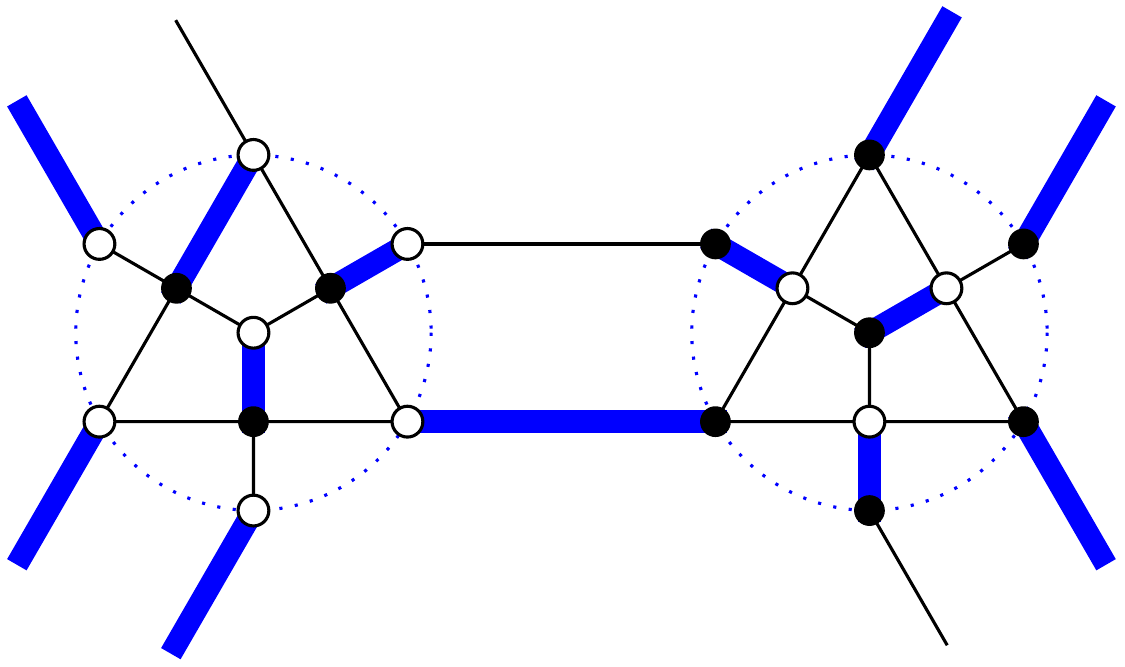}\hskip1cm\includegraphics[width=1.4in]{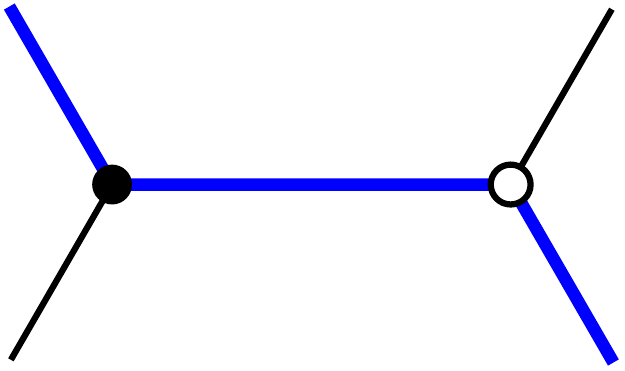}\\
        \vskip.5cm
        \includegraphics[width=2in]{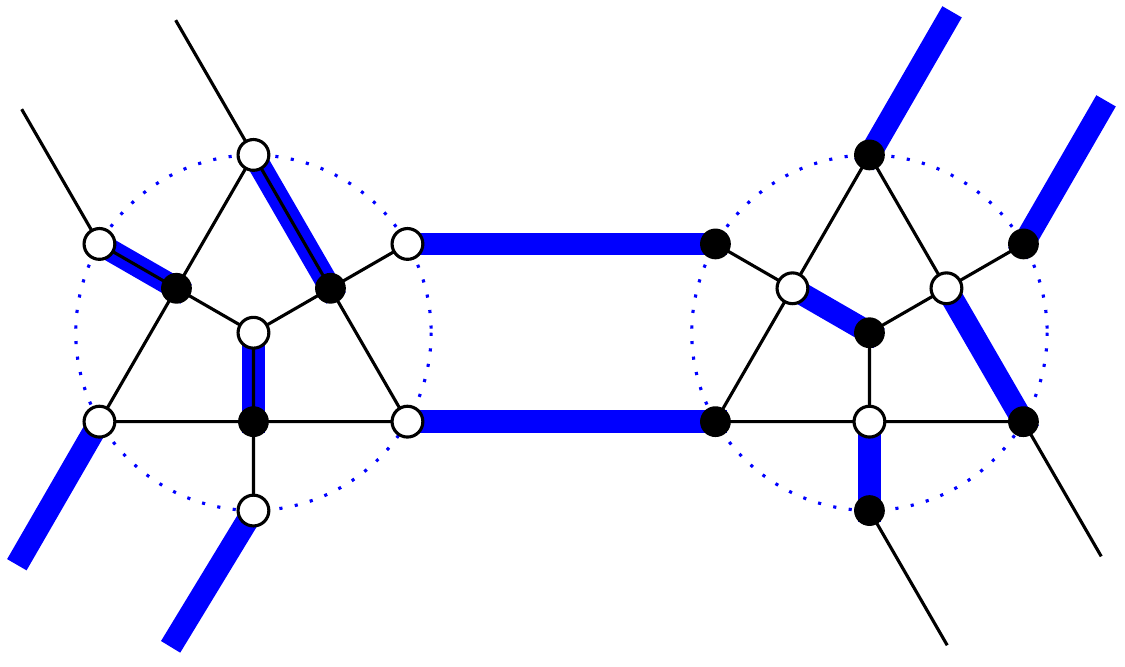}\hskip1cm\includegraphics[width=1.4in]{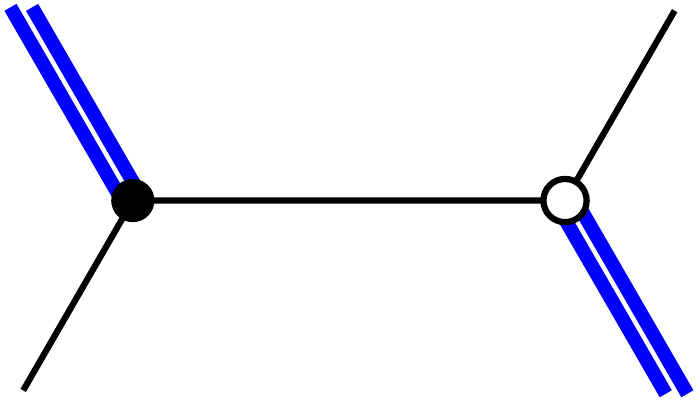}
    \end{center}
    \caption{\label{Gdimerex}
Mapping from dimer covers of $\hat\G$ (left) to $2$-multiwebs on $\G$ (right). Edges of $\G$ of multiplicity $m_e=0,1,2$ correspond to using respectively $2,1$ or $0$ edges of $\hat\G$ between the vertex gadgets. If $m_e=2$ or $m_e=0$ note that there is a unique (local) preimage. If $m_e=1$
there are multiple local preimages.}
\end{figure}

\begin{figure}
\begin{center} \includegraphics[width=1.6in]{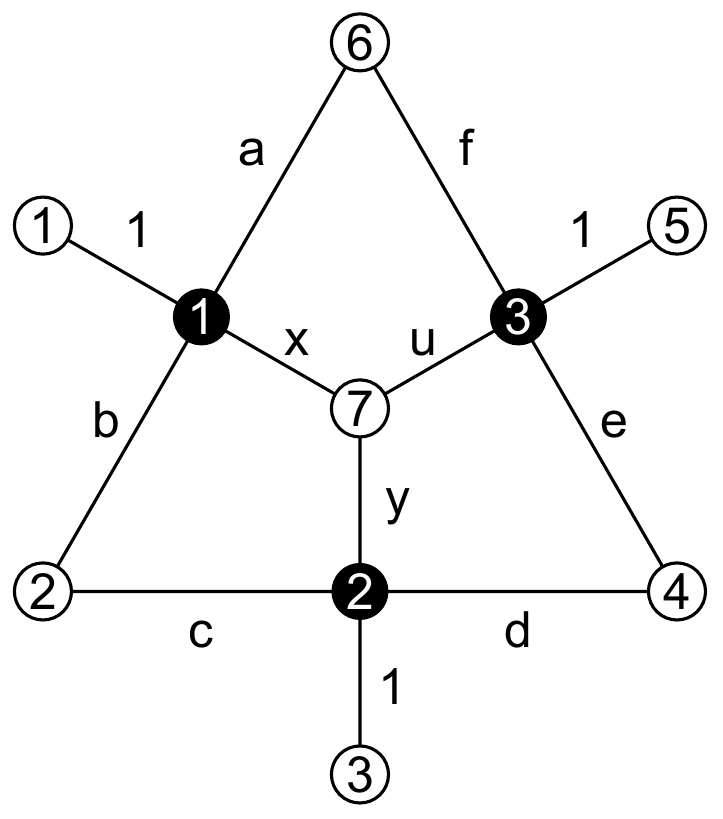} \hskip1cm\includegraphics[width=1.6in]{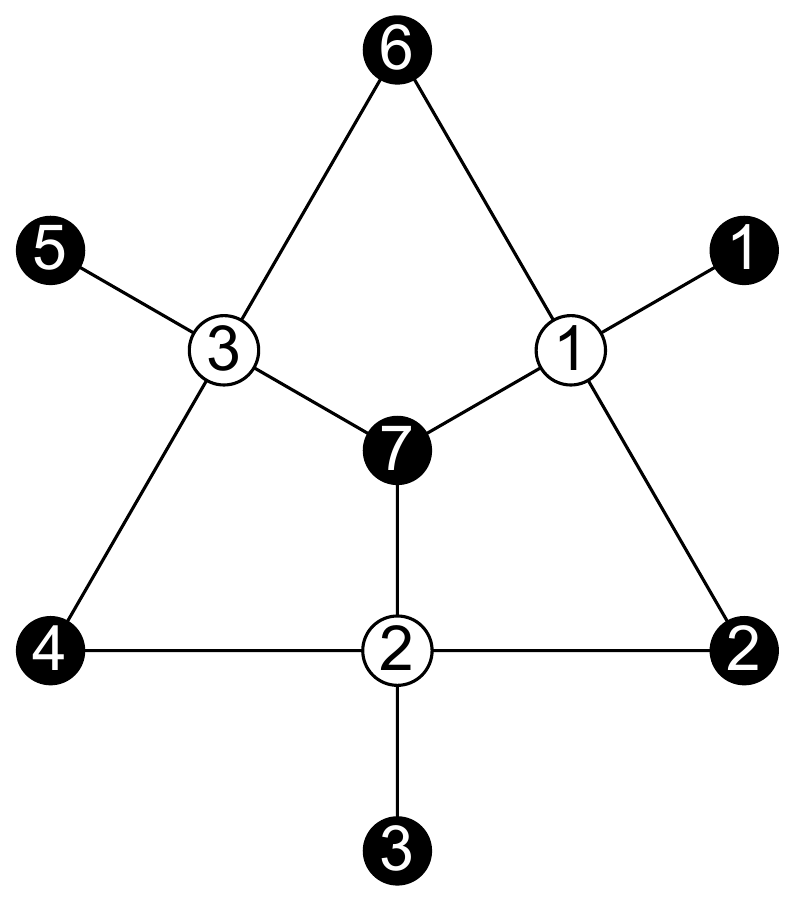}\end{center}
\caption{The ``gadget'' $H_{2,6}$ for a black vertex (left) and $H_{2,6}^*$ for a white vertex (right).
On the right, all edge weights are $1$. Note also that the index order of boundary vertices is counterclockwise for $H_{2,6}$, and clockwise for $H_{2,6}^*$.\label{H26}}
\end{figure}

Let $\G=(W\cup B,E)$ be a trivalent bipartite planar graph, and $\n\equiv 2$.
Let $\Phi=\{\phi_{bw}\}$ be a $\GL_2$-connection.
Let $\hat\G$ be the graph obtained from $\G$ by replacing each edge with two parallel edges and each vertex with a 
``gadget": the $10$-vertex graph $H_{2,6}$ or $H_{2,6}^*$ (depending on whether the vertex
is black or white) which parameterizes $\Grtp_{2,6}$, as shown in Figure 
\ref{hatG2}. We put edge weights as shown in Figure \ref{H26} on the gadgets $H_{2,6}$,
and edge weights $1$ on the gadgets $H_{2,6}^*$. 
The mapping from dimer
covers ($1$-multiwebs) on $\hat\G$ to $2$-multiwebs on $\G$ is illustrated in Figure \ref{Gdimerex}.

Note that $H_{2,6}$ is not exactly the graph one obtains from the construction presented at the end of Section \ref{sec:grassmannians}.
There are, however, local transformations of plabic graphs (see \cite{postnikov_06}) preserving the boundary measurements,
which transform this graph into $H_{2,6}$.

Let us first focus on the gadget $H_{2,6}$ of Figure \ref{H26}. With these edge weights, and a particular choice of Kasteleyn signs, the Kasteleyn matrix for the gadget is
\[ 
    K = \begin{pmatrix}
 1 & 0 & 0 \\
 b & c & 0 \\
 0 & 1 & 0 \\
 0 & d & e \\
 0 & 0 & 1 \\
 -a & 0 & f \\
 x & -y & u
    \end{pmatrix}.
\]

Here the last row corresponds to the white vertex internal to the gadget; the other white vertices communicate with the rest of the graph outside the gadget.
We can define a boundary response matrix by taking the Schur reduction of this $K$ onto its 
first six rows and first two columns. The reduced Kasteleyn matrix onto columns $1,2$ and rows $1,\dots,6$ is 
given (up to a scalar multiple of $z$ in the first column) by
\be\label{H26mat}
    K_{\mathrm{red}} = \begin{pmatrix}
        u & bu & 0 & -ex & -x & -(au+fx) \\
        0 & c & 1 & \frac{1}{u}(du+ey) & \frac{y}{u} & \frac{fy}{u}
\end{pmatrix}^t.
\ee
Choosing a different set of $2$ columns would result in an equivalent matrix: one with the same maximal minors, up to a global scale. 
The matrix $K_{\mathrm{red}}$ represents a point in the positive Grassmannian $\Grtp_{2,6}$ on condition that all weights $a,\dots,f,x,y,u$ are positive. 
Conversely any element of $\Grtp_{2,6}$ is equivalent (after column operations) to (\ref{H26mat}) with some choice of 
positive weights $a,\dots,f,x,y,u$. Letting these weights vary over $\R$ the above represents a generic element of $\Gr_{2,6}$. 

Likewise putting edge weights $1$ on the gadget $H_{2,6}^*$ leads to its boundary measurement
matrix 
\be\label{rowgrass}\begin{pmatrix}1&1&0&-1&-1&-2\\0&1&1&2&1&1\end{pmatrix},
\ee
the transpose of that for $H_{2,6}$ (with all variables set to $1$).
Finally, we need to find edge weights $a,b,\dots,u$ so that 
the matrix $\phi_{bw}$ on $\G$ corresponds to the $2\times 2$ matrix
obtained by multiplying the first two rows of the $H_{2,6}$ gadget at $b$ with the first two columns of the $H_{2,6}^*$ gadget at $w$: in this case
$$\phi_{bw} = \begin{pmatrix}u&0\\bu&c\end{pmatrix} \begin{pmatrix}1&1\\0&1\end{pmatrix}=\begin{pmatrix}u&u\\bu&bu+c\end{pmatrix}.$$
A similar expression holds for the other edges of $\G$.

Given matrices $\phi_{bw}$ on edges of $\G$, we therefore need to find edge weights on the gadget $H_{2,6}$ at each black vertex so that the appropriate rows of the Grassmannian element there, when multiplied by the corresponding columns of (\ref{rowgrass}), makes the $\phi_{bw}$.

We can guarantee this if we show that in fact for any three matrices $A,B,C\in GL_2(\R)$, the 
Grassmannian element $(A~B~C)^t$ can be realized with some choice of edge weights on $H_{2,6}$.
For a given generic Grassmannian element $(A~ B ~C)^t$, the edge weights on $H_{2,6}$ can be given explicitly in terms
of the matrices $A,B,C$ as follows. We have
\[ x = \Delta_{35}, \quad \quad \quad y = \Delta_{15}, \quad \quad \quad u = \Delta_{13} \]
\[ a = \frac{\Delta_{56}}{\Delta_{15}}, \quad \quad \quad c = \frac{\Delta_{12}}{\Delta_{13}}, \quad \quad \quad e = \frac{\Delta_{34}}{\Delta_{35}} \]
\[ b = \frac{\Delta_{23}}{\Delta_{13}}, \quad \quad \quad d = \frac{\Delta_{45}}{\Delta_{35}}, \quad \quad \quad f = \frac{\Delta_{16}}{\Delta_{15}} \]

These formulas hold for generic $A,B,C$. If some minors of $(A~B~C)^t$ are zero,
we can
take a limit of the corresponding generic cases, resulting in a simpler gadget with fewer
edges and vertices. More concretely, we find ourselves in a lower-dimensional positroid cell; Postnikov \cite{postnikov_06} shows how to assign a simpler gadget, depending on which minors are zero, which fully parameterizes such matrices. The analogs of the above formulas,
for the edge weights in terms of the minors,
for any positroid cell in $\Grnn_{k,n}$ were obtained in \cite{talaska}.

\subsection{Gadgets}

\begin {figure}[h!]
\centering
    \includegraphics[width=2in]{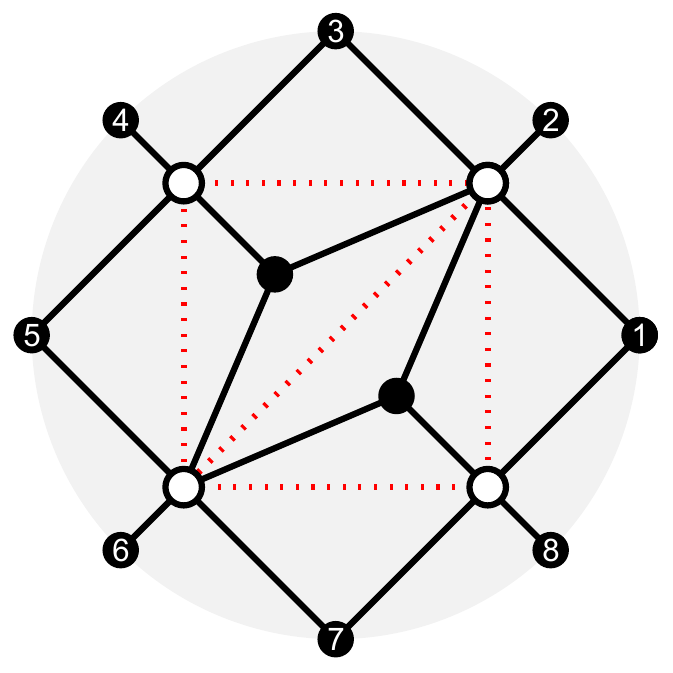}\hskip1cm\includegraphics[width=2in]{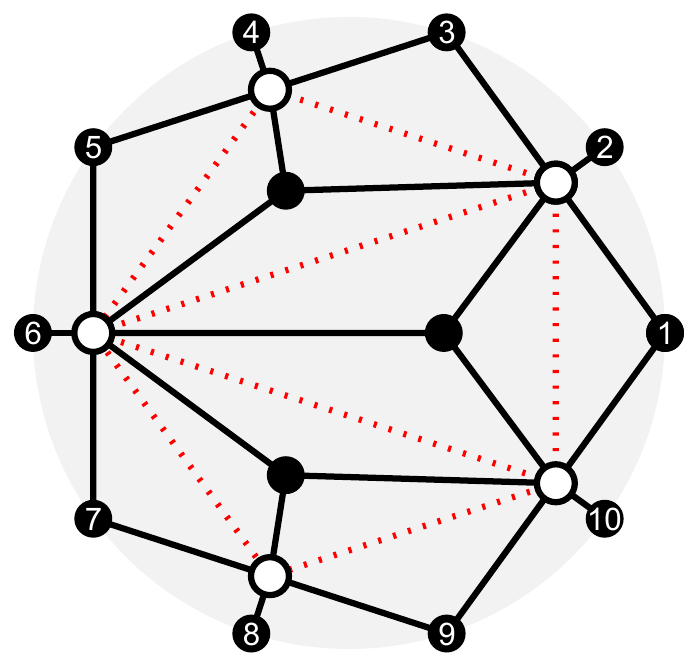}
\caption {Gadgets $H^*_{2,8}$ and $H^*_{2,10}$. These are constructed starting from a
triangulation of a square and pentagon, shown in dotted red. Gadgets for $H_{2,2n}^*$ can be 
constructed similarly from a $2n$-gon surrounding a triangulation of an $n$-gon;
equivalent gadgets are obtained by choosing other triangulations.}
\label {fig:gadgets}
\end {figure}

We define $H_{k,n}$ to be a reduced plabic graph representing the top-dimensional
cell in $\Grnn_{k,n}$ with white boundary vertices, and $H_{k,n}^*$ the same graph but with the colors reversed. We will refer
to these $H_{k,n}$ and $H_{k,n}^*$ graphs as ``gadgets''. Some examples can be seen in Figures \ref{H26} and \ref{fig:gadgets}.

We start with a graph $\G$ with vertex multiplicities $\n$ and a connection $\Phi$, and we will describe how to construct $\hat\G$ by
replacing each vertex with a gadget, and replacing each edge by multiple parallel edges connecting the gadgets. The rule is different at white vertices and black
vertices. Each edge $bw$ of $\G$ is replaced by $n_w$ parallel edges, where $n_w$ is the multiplicity at the white endpoint of the edge. Each white vertex $w$
of $\G$ with $d$ neighbors is replaced by a gadget $H^*_{n_w,dn_w}$; each black vertex $b$
is replaced by a gadget $H_{n_b,N}$ where $N=\sum_{w\sim b}n_w$. In this way for each edge $bw$ the gadget
at $b$ has an ordered subset of boundary vertices, of size $n_w$, corresponding
and connected via $n_w$ parallel edges, to the subset of size $n_w$ of the gadget at $w$.

We will assign all edges in gadgets $H_{n_w,dn_w}^*$ at white vertices of $\G$ a weight of $1$. 
Let $M_w$ be the reduced Kasteleyn matrix of the gadget at $w$; this is an $n_w\times dn_w$ matrix.
We write $M_w=[M_{b_1w},M_{b_2w},\dots,M_{b_dw}]$ its block decomposition according to the $d$ neighbors.

If an edge $bw$ of $\G$ has connection matrix $\phi_{bw}$, and if the corresponding block
of $M_w$ at $w$ is the matrix $M_{bw}$, then we assign the edge weights in the gadget at $b$ so that the appropriate block of its boundary measurement
matrix $M_b$ is the matrix $\phi_{bw}':=M_{bw}^{-1} \phi_{bw}$. That is, we write $\phi_{bw}=\phi_{bw}'M_{bw}$ where $M_{bw}$ is the 
fixed matrix coming from the gadget at $w$.

As in the basic example from the previous section, there is a mapping from dimer covers of the new graph $\hat\G$ onto the set of $\n$-multiwebs of 
the original graph $\G$ (pictured in Figure \ref{Gdimerex}). 
The image is the unique $\n$-multiweb on $\mathcal{G}$ whose multiplicities $m_e$ are the 
number of \emph{unoccupied} edges of the set of parallel edges of $\hat\G$ between $w$ and $b$.

\subsection {Gauge Transformations}

Let $\mathcal{G}$ be a bipartite graph with vertex labels $\n$ and a connection $\Phi = \{\phi_{bw}\}$. For a black vertex $b$ 
with $\n(b) = n_b$, and an invertible matrix $g \in \mathrm{GL}_{n_b}$, a \emph{gauge transformation} at $b$ changes the connection matrices
on all incident edges $bw$ by the rule $\phi_{bw} \mapsto \phi_{bw}g$. Similarly, for a white vertex $w$ with $\n(w) = n_w$ and some $g \in \mathrm{GL}_{n_w}$,
the corresponding gauge transformation alters weights of incident edges by $\phi_{bw} \mapsto g \phi_{bw}$. The trace of every multiweb on $\mathcal{G}$
is then multiplied by the scalar $\det(g)$.

Recall that for a $k$-valent white vertex $w$ with $\n(w) = n_w$, the corresponding gadget is $H^*_{n_w,kn_w}$, whose boundary measurement matrix
is a block matrix $K_w = (M_1, M_2, \dots, M_k)$. 
Gauge transformations at $w$ therefore correspond to
the left action of $\mathrm{GL}_{n_w}$ (i.e. $K_w \mapsto g K_w$). This agrees with the idea that the boundary measurement matrix should
be thought of as an element of the Grassmannian $\mathrm{Gr}_{n_w,kn_w}$.

Similarly, for a $k$-valent black vertex $b$ with $\n(b) = n_b$ and $N = \sum_{w \sim b} n_w$, the corresponding gadget is $H_{n_b,N}$, whose boundary
measurement matrix is the block matrix
\[ K_b = \begin{pmatrix} M_1^{-1} \phi_{b,w_1} \\ M_2^{-1} \phi_{b,w_2} \\ \vdots \\ M_k^{-1} \phi_{b,w_k} \end{pmatrix} \]
The gauge action at $b$ therefore corresponds to right multiplication (i.e. $K_b \mapsto K_b \cdot g$).

\subsection{Measure} \label{sec:measure}

For a planar graph $\G$ with vertex multiplicities $\n$ and connection $\Phi$, we constructed above a graph $\hat\G$,
the graph with gadgets, 
and a map $\Omega_1(\hat\G)\to\Omega_{\n}(\G)$. 

\begin{thm}\label{measurepres}The map $\Omega_1(\hat\G)\to\Omega_{\n}(\G)$ preserves the measure.
More precisely, the sum of weights of dimer covers of $\hat\G$ whose image is a given multiweb
$m$, is $\Tr(m)$. 
\end{thm}

\begin{proof}
Let $m\in\Omega_{\n}(\G)$. As before we split each edge $e$ of $\G$ into $m_e$ parallel
edges, so that  $m$ is simple, that is, all edge multiplicities are $1$ or $0$. 

It is convenient to put two new vertices on each edge of $\G$, both with multiplicity $n_w$,
so that original edge $bw$ with parallel transport $\phi_{bw}$ becomes $bw'b'w$, with parallel
transports $\phi'_{bw},I,$ and $M_{bw}$ on subedges $bw', w'b'$ and $b'w$ respectively. 
Since $\phi_{bw}=\phi'_{bw}M_{bw}$
this operation does not change $\Tr(m)$ for any multiweb $m$. 

Recall that in Section \ref{tracesection} we defined $\Tr(m)$ as the contraction of a certain tensor product of codeterminants. 
At a black vertex $b$, let $M_b$ be the element of $\Grtp_{n_b,N}$ there, written as an $N\times n_b$ matrix (that is,
as the transpose of the usual representation). 
Let $v_b$ be the codeterminant, for the web $m$, at $b$. Associated to the web $m$ is a 
submatrix $M_b(m)$ of $M_b$: the submatrix with rows corresponding to edges of $m$
(and remember that each edge of $m$ corresponds to multiple rows of $M_b$). We can apply $M_b(m)$ to the 
codeterminant $v_b$, which means applying, for each edge $bw'$ of $m$, the linear map $\phi'_{bw}$ to the entries of the tensor factor of $v_b$ for that edge. 
This results in a ``modified codeterminant" $M_b(m)v_b$. 

This modified codeterminant can be written as follows. Let $n=n_b$.
Starting from (\ref{codet}), apply $A=\phi'_{bw}$ to get
\begin{align*}M_b(m)v_b &= \sum_{\sigma\in S_n}(-1)^\sigma Au_{\sigma(1)}\otimes \dots \otimes Au_{\sigma(n)}\\
&= \sum_{\sigma\in S_n}(-1)^\sigma (\sum_{j_1} A_{j_1,\sigma(1)}u_{j_1})\otimes \dots\otimes(\sum_{j_k} A_{j_n,\sigma(n)}u_{j_n}).\end{align*}
The coefficient of $u_{i_1}\otimes\dots\otimes u_{i_n}$ in $M_b(m)v_b$ is thus
$$\sum_{\sigma\in S_k}(-1)^\sigma A_{i_1,\sigma(1)}\dots A_{i_n,\sigma(k)} = \det[(M_b)_{\{i_1,\dots,i_n\}}],$$
the maximal minor of $M_b$ corresponding to rows $\{i_1,\dots,i_n\}$.

Similarly at a white vertex $w$ we let $M_w$ be the element of $\Grtp_{n_w,dn_w}$ there, written
this time as an $n_w\times dn_w$ matrix. It has a submatrix $M_w(m)$ with columns corresponding to edges in $m$ (and each edge of $m$ corresponds to multiple columns). We apply 
$M_w(m)$ to the codeterminant $v_w$, by applying $(M_{bw})^t$ to the entries in the tensor factor for edge $b'w$
(we use the transpose because this is the action on the dual space). 
This results in a modified codeterminant $M_w(m)v_w$. We contract the tensor product of all the $M_b(m)v_b$ with
the tensor product of all the $M_w(m)v_w$'s to get the trace of $m$.

Now consider a dimer cover $\pi\in\Omega_1(\hat\G)$. For each edge $e$ of $\G$, it uses a subset $S(e,\pi)$ of edges of $\hat\G$ lying over the middle subedge of $e$. Consider the set $[\pi]\subset\Omega_1(\hat\G)$ of all dimer covers $\pi'$ using the same sets $S(e,\pi)$ for 
all edges $e$. This set is obtained by completing the dimer cover within each gadget. For a black vertex $b$ of $\G$,
suppose $\pi$ uses middle edges (that is, inputs) with indices $S_b = \{i_1,\dots\,i_n\}$ around the gadget at $b$.
The sum of weights of 
dimer covers completing the gadget
is the corresponding minor $I=\det[(M_b)_{\{i_1,\dots,i_n\}}]$. 
This is, as explained above, the coefficient of $u_{i_1}\otimes\dots\otimes u_{i_n}$ in the modified codeterminant
$M_b(m)v_b$. 

Likewise at a white vertex the sum of weights of dimer covers completing the gadget and using inputs 
$\{j_1,\dots,j_\ell\}$ is the coefficient of $f_{j_1}\otimes\dots\otimes f_{j_\ell}$ in the modified codeterminant
$M_w(m)v_w$. 

When contracting the modified codeterminants we get a nonzero contribution for every term in which
the indices on each edge of $\G$ match, in which case the contribution to the trace is the product of the 
coefficients, which is the product of the corresponding minors of the $M_v$. 

This completes the proof. 
\end{proof}

\section{$2$-multiwebs on honeycomb graphs}

Let $H$ be the bi-infinite honeycomb graph in $\R^2$, with 
constant vertex multiplicities $\n\equiv 2$ and with periodic 
$\text{GL}_2$-connection having matrices $A,B,C$ on the three directions of edges. 
Let $\H$ be a finite subgraph of $H$.
Under what conditions on $A,B,C$ do all $2$-multiwebs in $\H$ have positive trace?
We give a large class of examples here (see Theorem \ref{p2}).

\subsection{The parameter space} 

\begin {figure}[h!]
    \centering
    \begin {tikzpicture}
        \draw (-1,-1) -- (1,-1) -- (1,1) -- (-1,1) -- cycle;

        \draw (-0.25,-0.25) -- (0.25,0.25);
        \draw (0.25,0.25) -- (0,1);
        \draw (0,-1) -- (-0.25,-0.25);
        \draw (-0.25,-0.25) -- (-1,0);
        \draw (0.25,0.25) -- (1,0);

        \draw[blue] (-0.25,-0.25) -- (-0.25,0);
        \draw[blue] (0.25,0.25) -- (0,0.25);

        \draw[fill=white] (-0.25,-0.25) circle (0.06);
        \draw[fill=black] (0.25,0.25) circle (0.06);

        \draw (0,0) node[above left] {$I$};
        \draw (-0.1,-0.7) node[right] {$B$};
        \draw (-0.7,-0.1) node[below] {$A$};
    \end {tikzpicture}
    \caption {A simple bipartite graph on a torus}
    \label {fig:simple_torus_graph}
    \end {figure}
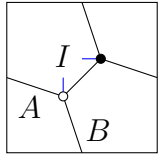

Consider the quotient $\G$ of $H$ by the lattice of its translational symmetries. This is the graph $\G$ on the torus with two vertices as pictured in Figure \ref{fig:simple_torus_graph}, with multiplicity $n=2$ at both vertices. 

Up to gauge transformations, we can assume the central edge has the identity connection, and the remaining edges have connection given by $2 \times 2$
matrices $A$ and $B$, as in Figure \ref{fig:simple_torus_graph}. 
For generic $A$ we can further apply a gauge transformation diagonalizing $A$. That is,
replace $A,B$ with $SAS^{-1}, SBS^{-1}$ where $SAS^{-1}$ is diagonal. Finally,
conjugating $B$
by any matrix commuting with $A$ (that is, any diagonal matrix) removes one more degree of freedom from the choice of $B$.
This leaves a $5$ dimensional space of matrix pairs $A,B$ which determine the system up to gauge. 

Associated to $\G$ is a polynomial 
$$P(z,w) = \det \tilde K = \det(I+zA+wB),$$
the \emph{characteristic polynomial}, see \cite{KOS}. We have
\be\label{cp}P(z,w) = 1+z^2\det A+w^2\det B+ z\,\tr A+w\,\tr B +zw \,\tr AB^*\ee
where $B^*=B^{-1}\det B$.

The characteristic polynomial describes the equivalent single-dimer system: see below.
It depends on a choice of basis for $H_1(\T^2,\R)$: $z$ and $w$ are the monodromies along the two basis cycles. 
Moreover it is only well-defined up to sign changes 
$z\mapsto \pm z$
and $w\mapsto \pm w$: these operations corresponds to different choices of Kasteleyn signs for the torus.

\subsection{Scalarization}

The scalarization $\hat{\mathcal{G}}$ is the one pictured in Figure \ref{hatG2} and Figure \ref{Gdimerex}.
This scalarization induces a scalarization of any subgraph $\H$ of the honeycomb $H$; 
on boundary vertices of $\H$ we ``cap off''
unused edges by removing them and their adjacent vertices from the corresponding gadgets.

The boundary measurement matrix at the right/white gadget (which has all its edge weights equal to 1) is
\[ M_w = \left( \begin{array}{cc|cc|cc} 1&0&-1&-3&-2&-3 \\ 0&1&1&2&1&1 \end{array} \right) = (X_1|X_2|X_3). \]
(This differs from the one of (\ref{rowgrass}) mentioned in Section \ref{sec:worked_example} by an $\mathrm{SL}_2$ change of basis, but this
does not change any of the $2 \times 2$ minors.) The boundary measurement matrix of the left/black gadget, using the edge weights as in Figure \ref{H26}, is
\[ 
    M_b = \begin{pmatrix}
        1 & 0 \\
        0 & 1 \\ \hline
        -bu & u \\
        -(cex+bdu+bey) & du+ey \\ \hline
        -(cx+by) & y \\
        -(acu+cfx+bfy) & fy
    \end{pmatrix}
    = \begin{pmatrix} Y_1 \\ \hline Y_2 \\ \hline Y_3 \end{pmatrix}
\]
Again, this differs from the one of (\ref{H26mat}) presented in the earlier section by an $\mathrm{SL}_2$ change of basis, and further using the assumption that $cu=1$.

The assumption that $\hat{\mathcal{G}}$ is the scalarization of $\mathcal{G}$ is equivalent to saying that $X_3Y_2 = A$ and $X_2Y_3 = B$.
This means the edge weights can be expressed as functions of the minors of $M$ (as given in Section \ref{sec:worked_example}), 
whose blocks are $Y_1 = \mathrm{Id}$, $Y_2 = X_3^{-1}A$ and $Y_3 = X_2^{-1}B$.

\begin {ex} \label{ex:all_ones}
    If all edge weights in both the white and black gadgets are set equal to $1$, then $M_b = M_w^\top$, and we get
    \[ A = \begin{pmatrix} 11 & -8 \\ -4 & 3 \end{pmatrix}, \quad \quad B = A^\top = \begin{pmatrix} 11 & -4 \\ -8 & 3 \end{pmatrix} \]
The characteristic polynomial (\ref{cp}) in this case is
$$P(z,w) = 1 + z^2 + w^2 + 14z+14w-14zw.$$
\end {ex}

\begin {ex}\label{trivexample}
    The trivial connection (when $A = B = \mathrm{Id}$) is realized in $\hat{\mathcal{G}}$ by the choice of weights
    (some of which are negative)
    \[ a = c = \frac{1}{3}, \quad b=d=e=f = -\frac{1}{3}, \quad x = -3, \quad y=u=3 \]
Its characteristic polynomial (\ref{cp}) in this case is
$$P(z,w) = 1 + z^2 + w^2 + 2z+2w+2zw=(1+z+w)^2.$$
\end {ex}

\subsection{Positive weights}

If the scalarization of $\G$ has positive edge weights, then all multiwebs on $\H$ will have positive 
traces, by Theorem \ref{measurepres}. This is a sufficient condition for positive traces (but not necessary; see example
\ref{trivexample} above).
We can determine when this happens as follows.

\begin{figure}
\begin{center}
\includegraphics[width=3in]{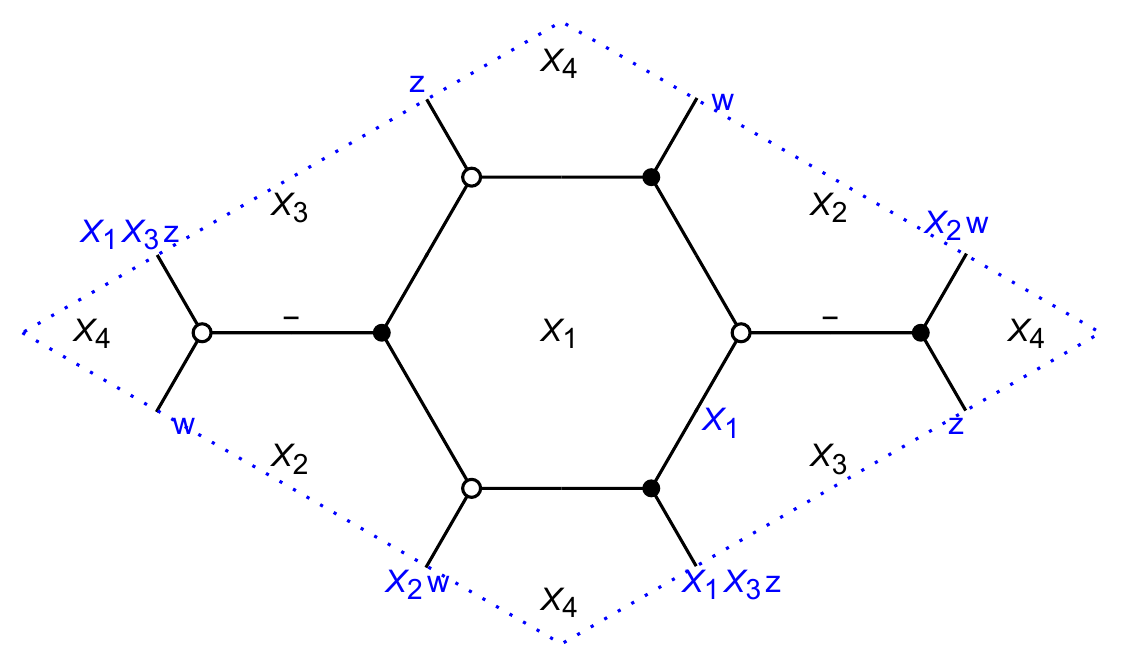}
\end{center}
\caption{Reduced graph corresponding to scalarization $\hat{\mathcal{G}}$ of the graph $\mathcal{G}$ from Figure \ref{fig:simple_torus_graph}.}
\label{2X2honeycomb}
\end{figure}

The graph of Figure \ref{hatG2} is not \emph{reduced}; by \cite{gk_13} it is equivalent, 
under local operations and gauge transformations, to a reduced graph
with the same characteristic polynomial, up to a multiplicative factor. 
The reduction of this graph is the one pictured in Figure \ref{2X2honeycomb},
which has $6$ parameters $X_1,X_2,X_3,X_4,z,w$, with the single relation $X_1X_2X_3X_4 = 1$. 
This reduced graph thus realizes the fact that the dimension of the space of parameters is indeed 5.

We can compute the reduction of the graph $\G$ directly by comparing the characteristic polynomials
of $\G$ and the graph in Figure \ref{2X2honeycomb}.
The partition function $\det\widetilde{K}$ must be equal, after rescaling $z,w$, and possibly changing their signs, to the characteristic
polynomial $P_{\mathrm{red}}(z,w)$ of the graph of  Figure \ref{2X2honeycomb}, that is, to
$$P_{\mathrm{red}}(z,w) = 1+(1+X_1X_3)z+X_1X_3z^2+(1+X_1X_2)w+X_1X_2w^2-(X_1+X_1X_2X_3)zw.$$

Comparing coefficients of $P(z,w)$ from (\ref{cp}) with $P_{\mathrm{red}}(\pm az,\pm bw)$ 
determines $A$ and $B$ up to gauge as functions of the five positive parameters $X_1,X_2,X_3,a,b$. 
We should allow for all four choices of signs.
For convenience let $\alpha=\pm a,\beta=\pm b$, so that $P_{\mathrm{red}}(\pm az,\pm bw)=P(\alpha z,\beta w)$.  
We find in particular the eigenvalues of $A$, $B$, and $AB^{-1}$ as functions of these parameters:
the eigenvalues of $A$ are $\alpha$ and $\alpha X_1X_3$, the eigenvalues of $B$ are $\beta$ and $\beta X_1X_2$, and the eigenvalues
of $AB^{-1}$ are $-\frac{\alpha}{\beta} X_2^{-1}$ and $-\frac{\alpha}{\beta} X_3$. This allows for the following observation.

\begin {thm}\label{p2}
  Let $\H$ be a finite balanced subgraph of the honeycomb graph $H$, with matrices $A,B$ and $I$ on the edges
    as in Figure \ref{fig:simple_torus_graph}, and periodic cilia as shown.
    Then the traces of all 2-multiwebs in $\H$ are positive if
    among the three matrices $A,B$ and $AB^{-1}$, either two have positive eigenvalues and one has negative eigenvalues,
    or all three have negative eigenvalues.
\end {thm}
 
\begin {proof}
Up to changing sign of $\alpha$ and/or $\beta$ (which corresponds to changing the sign of $A$ and/or $B$) we can suppose the eigenvalues of $A,B$ are positive and those of $AB^{-1}$ are negative. 
Then by the discussion above, we can choose positive values for
the five parameters $X_1$, $X_2$, $X_3$, $a$, and $b$, such that $\hat\H$ 
has the same dimer partition function as $\H$. By the measure-preserving property (see Section \ref{sec:measure}), the trace of
any 2-multiweb in $\H$ is equal to the weighted sum of dimer covers of $\hat{\H}$ with the same topological type as the multiweb. Since all weights
are positive, this means all multiweb traces are positive.
\end {proof}

There are other choices of matrices $A,B$ leading to positive traces, not covered by the theorem:
for example if $A,B$ are upper triangular with positive determinants.
We conjecture that these are the only two cases:
\begin{conj}
Traces of all $2$-webs in $H$ are positive if and only if, up to gauge equivalence, one of the following conditions holds:
\begin{enumerate}
\item $A$ and $B$ are upper triangular with positive determinants, or
\item Among the three matrices $A,B$ and $AB^{-1}$, either two have positive eigenvalues and one has negative eigenvalues, or all three have negative eigenvalues.
\end{enumerate}
\end{conj}

\old{
Conversely, suppose all multiweb traces in $G$ are positive. In this small example, there are only six multiwebs, whose traces are $1$, $\tr(A)$, $\tr(B)$,
    $\det(A)$, $\det(B)$, and $\det(B)\tr(AB^{-1})$. Knowing $\det(A) > 0$ means the eigenvalues of $A$ have the same sign, and knowing $\tr(A) > 0$ means they
    must both be positive (and similarly for $B$). Since $\det(AB^{-1}) = \frac{\det(A)}{\det(B)} > 0$, we know its eigenvalues have the same sign. Finally, the
    assumption $\tr(AB^{-1}) < 0$ implies both eigenvalues must be negative. 
}

\bigskip

\begin {ex}
    In the case from Example \ref{ex:all_ones}, where all edge weights are set equal to 1, the eigenvalues of $A$ and $B$ are both $7 \pm 4\sqrt{3}$,
    which are positive, and the eigenvalues of $AB^{-1}$ are $-7 \pm 4\sqrt{3}$, which are both negative.
    This corresponds to the network in Figure \ref{2X2honeycomb} with the weights
    \[ a = X_1 = X_2 = X_4 = 7 + 4\sqrt{3} \approx 13.928, \quad b = 7 - 4\sqrt{3} \approx 0.0718 \]
    \[ X_3 = (7-4\sqrt{3})^3 = 1351 - 780\sqrt{3} \approx 0.00037 \]
\end {ex}

\bigskip

\begin {ex}
    Similar to the previous example, we could instead set all the parameters in Figure \ref{2X2honeycomb} equal to 1, so that
    $X_1 = \cdots = X_4 = a = b = 1$. In this case the representative matrices are
    \[ A = \begin{pmatrix} 1 & 0 \\ 4 & 1 \end{pmatrix}, \quad \quad B = \begin{pmatrix} 1&1\\0&1 \end{pmatrix}, \quad \quad AB^{-1} = \begin{pmatrix} 1&-1\\4&-3\end{pmatrix} \]
    It is easy to see that $A$ and $B$ have eigenvalues $1,1$, and $AB^{-1}$ has eigenvalues $-1,-1$.
    The fact that $\tr(A) = 2$ reflects the fact that there are exactly 2 perfect matchings of the graph in Figure \ref{2X2honeycomb}
    with 1 edge crossing the $z$-boundary (and similarly for $B$ and $AB^{-1}$). The fact that $\det(A) = 1$ reflects the fact
    that there is one perfect matching using two edges which cross the $z$-boundary (and similarly for $B$).
\end {ex}

\bibliographystyle{alpha}
\bibliography{higher_networks}

\end {document}